\documentclass[leqno,12pt]{article}  
\usepackage{amsmath}
\usepackage{amssymb}

\usepackage[german,english]{babel}
\usepackage{amsthm}

\title{Integral structures in the $p$-adic holomorphic discrete series}

\author{\textsc{Elmar Grosse-Kl\"onne}}
\date{}

\theoremstyle{plain} 
\newtheorem{satz}{Theorem}[section]  
\newtheorem{kor}[satz]{Corollary}  
\newtheorem{lem}[satz]{Lemma}  
\newtheorem{pro}[satz]{Proposition}  
 
\newcommand{\spec}{\mbox{\rm Spec}}  
\newcommand{\proj}{\mbox{\rm Proj}}
\newcommand{\spf}{\mbox{\rm Spf}}  

\newcommand{\kara}{\mbox{\rm char}}  

\newcommand{\sym}{\mbox{\rm Sym}}

\newcommand{\dlog}{\mbox{\rm dlog}}

\newcommand{\diag}{\mbox{\rm diag}}

\theoremstyle{remark}

\theoremstyle{definition}

\DeclareMathOperator{\Hom}{Hom}

\newcommand{\0}{\ensuremath{\overrightarrow{0}}}

\begin{document}
\maketitle
\footnote[0]{I wish to thank Peter Schneider and Jeremy Teitelbaum for generously sharing their ideas on how in their work \cite{schtei}, integral structures could possibly be taken into account. I thank the referee for hints and comments on the text.}

\begin{abstract}

For a local non-Archimedean field $K$ we construct ${\rm GL}_{d+1}(K)$-equivariant coherent sheaves ${\mathcal V}_{{\mathcal O}_K}$ on the formal ${\mathcal O}_K$-scheme ${\mathfrak X}$ underlying the symmetric space $X$ over $K$ of dimension $d$. These ${\mathcal V}_{{\mathcal O}_K}$ are ${\mathcal O}_K$-lattices in (the sheaf version of) the holomorphic discrete series representations (in $K$-vector spaces) of ${\rm GL}_{d+1}(K)$ as defined by P. Schneider \cite{schn}. We prove that the cohomology $H^t({\mathfrak X},{\mathcal V}_{{\mathcal O}_K})$ vanishes for $t>0$, for ${\mathcal V}_{{\mathcal O}_K}$ in a certain subclass. The proof is related to the other main topic of this paper: over a finite field $k$, the study of the cohomology of vector bundles on the natural normal crossings compactification $Y$ of the Deligne-Lusztig variety $Y^0$ for ${\rm GL}_{d+1}/k$ (so $Y^0$ is the open subscheme of ${\mathbb P}_k^d$ obtained by deleting all its $k$-rational hyperplanes).    

\end{abstract}

%


\begin{center} {\bf Introduction} 
\end{center}

Let $K$ be a non-Archimedean locally compact field with residue field $k$ of characteristic $p$. Let $d\in\mathbb{N}$ and let $X$ be the Drinfel'd symmetric space over $K$ of dimension $d$: the $K$-rigid space which is the complement in $\mathbb{P}_K^d$ of all $K$-rational hyperplanes. The group $G={\rm GL}_{d+1}(K)$ acts on $X$, and it is expected that the cohomology of $G$-equivariant sheaves on $X$ affords wide classes of interesting $G$-representations in infinite dimensional vector spaces. For example, by now we know that the ${\ell}$-adic cohomology ($\ell\ne p$) of certain \'{e}tale coverings of $X$ contains all the {\it smooth} discrete series representations of $G$ (in characteristic zero). A very different class of $G$-representations in infinite dimensional $K$-vector spaces is obtained by taking the global sections of $G$-equivariant vector bundles on $X$. The study of these has been initiated by Morita in the case $d=1$. In that case the relevant vector bundles are automorphic line bundles on $X$ which are classified by their weight (an integer, or equivalently: an irredicible $K$-rational representation of ${\rm GL}_1$). Generalizing to any $d$ Schneider \cite{schn} assigns to an irredicible $K$-rational representation $V$ of ${\rm GL}_d$ a $G$-equivariant vector bundle ${\mathcal V}$ on $X$ (in \cite{schn} he in fact only considers the action by ${\rm SL}_{d+1}(K)$; here we will consider a suitable extension to $G$ which however depends on the choice of $\widehat{\pi}\in \widehat{K}$, see the text). The resulting $G$-representations ${\mathcal V}(X)$, which he called the "holomorphic discrete series representations", are at present very poorly understood, at least if $d>1$. Already the seemingly innocent case ${\mathcal V}=\omega_X$, the line bundle of $d$-forms on $X$, turned out to be fairly intricate and required a host of original techniques (Schneider and Teitelbaum \cite{ast}). 

It is natural to ask for integral structures inside ${\mathcal V}$. Let ${\mathfrak X}$ be the natural $G$-equivariant strictly semistable formal ${\mathcal O}_K$-scheme with generic fibre $X$ constructed in \cite{mus} and consider ${\mathcal V}$ as a $G$-equivariant sheaf on ${\mathfrak X}$ via the specialization map $X\to{\mathfrak X}$. Let $\widehat{K}/K$ be a totally ramified extension of degree $d+1$. In the case $d=1$ we constructed in \cite{mathan} a $G$-stable ${\mathcal O}_{\mathfrak X}\otimes_{{\mathcal O}_K}{\mathcal O}_{\widehat{K}}$-coherent subsheaf ${\mathcal V}_{{\mathcal O}_{\widehat{K}}}$ inside ${\mathcal V}\otimes_K{\widehat{K}}$, generalizing a previous construction of Teitelbaum from the case of even weight \cite{teit}. We completely determined the cohomology $H^*({\mathfrak X},{\mathcal V}_{{\mathcal O}_{\widehat{K}}})$ and obtained (if $\kara(K)=0$) applications to the $\Gamma$-group cohomology $H^*(\Gamma,{\mathcal V}(X))$ of the above $G$-representations ${\mathcal V}(X)$ over $K$, for cocompact discrete subgroups $\Gamma\subset{\rm SL}_2(K)$. Let us also mention that Breuil uses (among many other tools) these ${\mathcal V}_{{\mathcal O}_{\widehat{K}}}$ to construct certain Banach space representations of ${\rm GL}_2({\mathbb Q}_p)$ which he expects to occur in a hoped for $p$-adic continuous Langlands correspondence. 

In the present paper we construct integral structures ${\mathcal V}_{{\mathcal O}_{\widehat{K}}}$ inside ${\mathcal V}\otimes_K{\widehat{K}}$ for arbitrary $d$. Roughly we proceed as follows. We have ${\mathcal V}\otimes_K{\widehat{K}}=V\otimes_{{\mathcal O}_K}{\mathcal O}_{\mathfrak X}\otimes_{{\mathcal O}_K}{{\mathcal O}}_{\widehat{K}}$ as an ${\mathcal O}_{\mathfrak X}\otimes_{{\mathcal O}_K}\widehat{K}$-module. The definition of the action of $G$ is based on an embedding ${\rm GL}_d\to{\rm GL}_{d+1}$ which restricts to an embedding $T_1\to T$ of the respective tori of diagonal matrices. We fix a ${\rm GL}_d({\mathcal O}_K)$-stable ${\mathcal O}_K$-lattice $V_0$ inside $V$. It decomposes as $V_0=\oplus_{\mu}V_{0,\mu}$ with the sum running over the weights $\mu$ of $V$ with respect to $T_1$. We choose ${\mathfrak Y}$, an open formal subscheme of ${\mathfrak X}$ such that the set of irreducible components of the reduction  ${\mathfrak Y}\otimes k$ is an orbit for the action of $T$ on the set of irreducible components of ${\mathfrak X}\otimes k$. Restricted to ${\mathfrak Y}$ we define ${\mathcal V}_{{\mathcal O}_{\widehat{K}}}|_{{\mathfrak Y}}=\oplus_{\mu}({\mathcal V}_{{\mathcal O}_{\widehat{K}}}|_{{\mathfrak Y}})_{\mu}$ where $({\mathcal V}_{{\mathcal O}_{\widehat{K}}}|_{{\mathfrak Y}})_{\mu}$ is a ${\mathcal O}_{\mathfrak X}\otimes_{{\mathcal O}_K}{\mathcal O}_{\widehat{K}}$-coherent lattice inside $V_{0,\mu}\otimes_{{\mathcal O}_K}{\mathcal O}_{\mathfrak X}\otimes_{{\mathcal O}_K}{\widehat{K}}$ whose position relative to the constant ${\mathcal O}_{\mathfrak X}\otimes_{{\mathcal O}_K}{\mathcal O}_{\widehat{K}}$-lattice $V_{0,\mu}\otimes_{{\mathcal O}_K}{\mathcal O}_{\mathfrak X}\otimes_{{\mathcal O}_K}{\mathcal O}_{\widehat{K}}$ is given by $\mu$. Then we prove that there exists a unique $G$-stable extension ${\mathcal V}_{{\mathcal O}_{\widehat{K}}}$ of ${\mathcal V}_{{\mathcal O}_{\widehat{K}}}|_{{\mathfrak Y}}$ to all of ${\mathfrak X}$.

We begin to compute the cohomology $H^*({\mathfrak X},{\mathcal V}_{{\mathcal O}_{\widehat{K}}})$, but the results we obtain here are by no means complete. However, we obtain a clean result for $V$ which are "strongly dominant": if we identify as usual the weights $\mu$ of $V$ with vectors $(a_1,\ldots,a_d)\in\mathbb{Z}^d$, then we require $\sum_{i\ne j}a_i\le da_j$ for all $1\le j\le d$, for all weights. In particular, $V$ is dominant in the sense that $0\le a_i$ for all $1\le i\le d$, for all weights. \\

{\bf Theorem \ref{globnull}:} {\it Suppose that $V$ is strongly dominant. Then} \begin{gather}H^t(\mathfrak{X},{\mathcal V}_{{\mathcal O}_{\widehat{K}}})=0\quad\quad(t>0),\tag{$i$}\\H^t(\mathfrak{X},{\mathcal V}_{{\mathcal O}_{\widehat{K}}}\otimes_{{\mathcal O}_{\widehat{K}}}k)=0\quad\quad(t>0),\tag{$ii$}\\H^0(\mathfrak{X},{\mathcal V}_{{\mathcal O}_{\widehat{K}}})\otimes_{{\mathcal O}_{\widehat{K}}}k=H^0({\mathfrak{X}},{\mathcal V}_{{\mathcal O}_{\widehat{K}}}\otimes_{{\mathcal O}_{\widehat{K}}}k).\tag{$iii$}\end{gather}

If $d=1$, strong dominance is equivalent with (usual) dominance and our results in \cite{mathan} showed that for nontrivial $V$, dominance is equivalent with the validity of the vanishing assertions from Theorem \ref{globnull}. If however $d>1$, (usual) dominance is not enough to guarantee the vanishing assertions from Theorem \ref{globnull}; a first counterexample is the case where $d=2$ and $V$ has highest weight corresponding to the vector $(8,3)$.  

Examples of ${\mathcal V}_{{\mathcal O}_{\widehat{K}}}$'s arising from strongly dominant representations $V$ are the terms $\Omega^s_{\mathfrak X}$ of the relative logarithmic de Rham complex of ${\mathfrak X}\to\spf({\mathcal O}_K)$ (with respect to canonical log structures). Of course these are defined even over ${\mathcal O}_K$, for them the extension to ${\mathcal O}_{\widehat{K}}$ is unnecessary. Thus Theorem \ref{globnull} implies $H^t(\mathfrak{X},\Omega^s_{\mathfrak X})=0$ for all $t>0$ and $H^0(\mathfrak{X},\Omega^s_{\mathfrak X})\otimes k=H^0(\mathfrak{X},\Omega^i_{\mathfrak X}\otimes k)$ (all $s$). Statements $(i)$ and $(iii)$ of Theorem \ref{globnull} follow from statement $(ii)$. The proof of $(ii)$ is reduced (using the main result from \cite{acy}) to the vanishing of $H^t(Z,{\mathcal V}_{{\mathcal O}_{\widehat{K}}}\otimes_{{\mathcal O}_{\widehat{\mathfrak X}}}{\mathcal O}_Z)$ for $t>0$, for all irreducible components $Z$ of ${\mathfrak X}\otimes k$, where we write ${{\mathcal O}_{\widehat{\mathfrak X}}}={{\mathcal O}_{{\mathfrak X}}}\otimes_{{\mathcal O}_{{K}}}{\mathcal O}_{\widehat{K}}$.

A typical irreducible component $Y=Z$ of ${\mathfrak X}\otimes k$ is isomorphic to the natural compactification, with normal crossings divisors at infinity, of the complement $Y^0$ in ${\mathbb P}^d_k$ of all $k$-rational hyperplanes; explicitly, $Y$ is the successive blowing up of ${\mathbb P}^d_k$ in all its $k$-rational linear subvarieties. ${\rm GL}_{d+1}(k)$ acts on $Y$, and the study of the cohomology of ${\rm GL}_{d+1}(k)$-equivariant vector bundles on $Y$ (like ${\mathcal V}_{{\mathcal O}_{\widehat{K}}}\otimes_{{\mathcal O}_{\widehat{\mathfrak X}}}{\mathcal O}_Y$) is the analog over $k$ of the program over $K$ described above. Of course this study should be of interest in its own rights.

Here we begin (section \ref{cohmodp}) by establishing vanishing theorems for the cohomology of certain line bundles associated with divisors on $Y$ which are stable under the subgroup of unipotent upper triangular matrices in ${\rm GL}_{d+1}(k)$. Given the cohomology of line bundles on projective space this turned out to be an essentially combinatorical matter (the underlying object is the building of ${\rm PGL}_{d+1}/k$); in fact we enjoyed this computation. We obtain $H^t(Y,{\mathcal V}_{{\mathcal O}_{\widehat{K}}}\otimes_{{\mathcal O}_{\widehat{\mathfrak X}}}{\mathcal O}_Y)=0$ for $t>0$ as desired (for strongly dominant $V$). 

In principle our vanishing theorems for first cohomology groups allow the determination of the ${\rm GL}_{d+1}(k)$-representation on the finite dimensional $k$-vector space $H^0(Y,{\mathcal V}_{{\mathcal O}_{\widehat{K}}}\otimes_{{\mathcal O}_{\widehat{\mathfrak X}}}{\mathcal O}_Y)$, and hence of the $G$-representation on the infinite dimensional $k$-vector space $$H^0(\mathfrak{X},{\mathcal V}_{{\mathcal O}_{\widehat{K}}})\otimes_{{\mathcal O}_{\widehat{K}}}k.$$Note that  $H^0(\mathfrak{X},{\mathcal V}_{{\mathcal O}_{\widehat{K}}})$ is a $G$-stable ${{\mathcal O}_{\widehat{K}}}$-submodule of the holomorphic discrete series representation ${\mathcal V}(X)\otimes_K{\widehat{K}}$. We start this discussion in section \ref{secdera} by analyzing the ${\rm GL}_{d+1}(k)$-representations $H^0(Y,{\mathcal V}_{{\mathcal O}_{\widehat{K}}}\otimes_{{\mathcal O}_{\widehat{\mathfrak X}}}{\mathcal O}_Y)$ with ${\mathcal V}_{{\mathcal O}_{\widehat{K}}}=\Omega^s_{\mathfrak X}\otimes_{{\mathcal O}_K}{\mathcal O}_{\widehat{K}}$ for some $s$. If $\Omega_Y^{\bullet}$ denotes the de Rham complex on $Y$ with allowed logarithmic poles along $Y-Y^0$ then $\Omega_Y^{\bullet}\cong (\Omega^s_{\mathfrak X}\otimes_{{\mathcal O}_K}{\mathcal O}_{\widehat{K}})\otimes_{{\mathcal O}_{\widehat{\mathfrak X}}}{\mathcal O}_Y=\Omega^s_{\mathfrak X}\otimes_{{\mathcal O}_{{\mathfrak X}}}{\mathcal O}_Y$, thus we are studying the ${\rm GL}_{d+1}(k)$-representations $H^0(Y,\Omega_Y^s)$. We describe an explicit $k$-basis of $H^0(Y,\Omega_Y^s)$ consisting of logarithmic differential forms and show that, as a ${\rm GL}_{d+1}(k)$-representation, it is a generalized Steinberg representation. As a corollary of our explicit computations we obtain the irreducibility of these generalized Steinberg representations. Moreover we derive that the log crystalline cohomology $H^s_{crys}(Y/W(k))$ is a representation of ${\rm GL}_{d+1}(k)$ on a {\it free} finite $W(k)$-module, and that $H^0(Y,\Omega_Y^s)$ is its reduction modulo $p$. Note that $H^s_{crys}(Y/W(k))\otimes{\mathbb Q}=H_{rig}^s(Y^0)$, the rigid cohomology with constant coefficients of the Deligne-Lusztig variety $Y^0$ for ${\rm GL}_{d+1}/k$. 

Another application of our vanishing theorems can be found in \cite{latt}: they enable us to compute the cohomology of sheaves of bounded logarithmic differential forms on ${\mathfrak X}$, with coefficients in certain algebraic ${\rm GL}_{d+1}(K)$-representations. This leads to the proof of certain (previously unknown) cases of a conjecture of Schneider, formulated in \cite{schn}, concerning Hodge decompositions of the de Rham cohomology (with coefficients) of projective $K$-varieties uniformized by $X$.\\
   
{\it Notations:} $K$ denotes a non-Archimedean locally compact field and $K_a$ its algebraic closure, ${\cal O}_K$ its ring of integers, $\pi\in{\cal O}_K$ a fixed prime element and $k$ the residue field with $q$ elements, $q\in p^{\mathbb{N}}$. We denote by $\omega:K^{\times}_a\to\mathbb{Q}$ the extension of the discrete valuation $\omega:K^{\times}\to\mathbb{Z}$ normalized by $\omega(\pi)=1$. We fix $\widehat{\pi}\in K_a$ with $\widehat{\pi}^{d+1}=\pi$ and set $\widehat{K}=K(\widehat{\pi})$. 
 
We fix $d\in \mathbb{N}$ and enumerate the rows and columns of ${\rm GL}\sb {d+1}$-elements by $0,\ldots,d$. We let ${{{U}}}\subset{\rm GL}\sb {d+1}$ denote the subgroup of unipotent upper triangular matrices,$${{{U}}}=\{(a_{ij})_{0\le i,j\le d}\in{\rm GL}\sb {d+1}\quad|\quad\,a_{ii}=1\,\mbox{for all}\,i,\quad \,a_{ij}=0\,\mbox{if}\,i>j\}.$$We set $G={\rm GL}\sb {d+1}(K)$. For $r\in\mathbb{R}$ we define $\lfloor r\rfloor, \lceil r\rceil\in\mathbb{Z}$ by requiring $\lfloor r\rfloor\le r<\lfloor r\rfloor+1$ and $\lceil r\rceil-1<r\le\lceil r\rceil$. For a divisor $D$ on a smooth connected $k$-scheme $S$ we denote by ${\mathcal L}_S(D)$ the associated line bundle on $S$, endowed with its canonical embedding into the constant sheaf generated by the function field of $S$.

\section{Line bundle cohomology of rational varieties}
\label{cohmodp}

The (right) action of ${\rm GL}_{d+1}(k)={\rm GL}(k^{d+1})$ on $(k^{d+1})^*=\Hom_k(k^{d+1},k)$ defines a (left) action of ${\rm GL}_{d+1}(k)$ on the affine $k$-scheme associated with $(k^{d+1})^*$, and this action passes to a (left) action of ${\rm GL}_{d+1}(k)$ on the projective space $$Y_0={\mathbb P}((k^{d+1})^*)\cong {\mathbb P}_k^{d}.$$For $0\le j\le d-1$ let ${\mathcal V}_0^j$ be the set of all $k$-rational linear subvarieties $Z$ of $Y_0$ with $\dim(Z)=j$, and let ${\mathcal V}_0=\bigcup_{j=0}^{d-1}{\mathcal V}_0^j$. The sequence of projective $k$-varieties$$Y=Y_{d-1}{\longrightarrow}Y_{d-2}{\longrightarrow}\ldots{\longrightarrow}Y_0$$is defined inductively by letting $Y_{j+1}\to Y_j$ be the blowing up of $Y_j$ in the strict transforms (in $Y_j$) of all $Z\in {\mathcal V}_0^j$. The action of ${\rm GL}_{d+1}(k)$ on $Y_0$ naturally lifts to an action of ${\rm GL}\sb {d+1}(k)$ on $Y$. Let $\Xi_0,\ldots,\Xi_{d}$ be the standard projective coordinate functions on $Y_0$ and hence on $Y$ corresponding to the canonical basis of $(k^{d+1})^*$. Hence $Y_0=\proj(k[\Xi_i;\,i\in\Upsilon])$ with $$\Upsilon=\{0,\ldots,d\}.$$ For $\emptyset\ne\tau\subset\Upsilon$ let $V_{\tau,0}$ be the reduced closed subscheme of $Y_0$ which is the common zero set of $\{\Xi_i\}_{i\in\tau}$. Let $V_{\tau}$ be the closed subscheme of $Y$ which is the strict transform of $V_{\tau,0}$ under $Y\to Y_0$. These $V_{\tau}$ are particular elements of the following set of divisors on $Y$:$${\mathcal V}=\mbox{the set of all strict transforms in}\,\,Y\,\,\mbox{of elements of}\,\,{\mathcal V}_0.$$ For $\tau\subset\Upsilon$ let $$\tau^c=\Upsilon-\tau,$$$${{{U}}}_{\tau}=\{(a_{ij})_{0\le i,j\le d}\in{{{U}}}(k)\quad|\quad a_{ij}=0\,\mbox{if}\,i\ne j\,\mbox{and}\,[j\in\tau^c\,\mbox{or}\,\{i,j\}\subset \tau]\}.$$Let $${{{\mathcal{Y}}}}=\{\tau\quad|\quad\emptyset\ne\tau\subsetneq\Upsilon\},$$$${\mathcal N}=\{(\tau,u)\quad|\quad\tau\in{{{\mathcal{Y}}}}, u\in{{{U}}}_{\tau}\}.$$We have bijections$${\mathcal N}\cong{\mathcal V},\quad\quad (\tau,u)\mapsto u.V_{{\tau}},$$$${{{\mathcal{Y}}}}\cong(\mbox{the set of orbits of}\,\,{{{U}}}(k)\,\,\mbox{acting on}\,\,{\mathcal V}),\quad\tau\mapsto\{u.V_{{\tau}}|\,u\in{{{U}}}_{\tau}\}.$$

{\it Remark:} We also have a bijection between ${\mathcal N}$ and the set of vertices of the building associated to ${\rm PGL}\sb {d+1}/k$. The subset $\{(\tau,1)|\,\,\tau\in\mathcal{Y}\}$ of ${\mathcal N}$ is then the set of vertices in a standard apartment.\\ 

For an element $\overline{a}=(\overline{a}_1,\ldots,\overline{a}_d)$ of $\mathbb{Z}^d$ let $\overline{a}_0=-\sum_{j=1}^d\overline{a}_j$ and for $\sigma\in{{{\mathcal{Y}}}}$ let$$b_{\sigma}(\overline{a})=-\sum_{j\in\sigma}\overline{a}_j.$$Given two more elements ${n}=(n_1,\ldots,n_d)$ and ${m}=(m_1,\ldots,m_d)$ of $\mathbb{Z}^d$ we define the divisor\begin{gather}D(\overline{a},n,m)=\sum_{\sigma\in{\mathcal Y}\atop 0\notin\sigma}(m_{|\sigma|}+b_{\sigma}(\overline{a}))\sum_{u\in{{{U}}}_{\sigma}}u.V_{\sigma}+\sum_{0\in\sigma\in{\mathcal Y}}(n_{|\sigma|}+b_{\sigma}(\overline{a}))\sum_{u\in{{{U}}}_{\sigma}}u.V_{\sigma}\label{divdef}\end{gather}on $Y$. The purpose of this section is to prove vanishing theorems for the cohomology of line bundles on $Y$ of the type ${\mathcal L}_Y(D(\overline{a},n,m))$ for suitable $\overline{a},n,m$.\\

We need relative (or restricted) analogs of the above definitions. For a subset $\sigma\subset\Upsilon$ we define the sequence of projective $k$-varieties$$Y^{\sigma}=Y^{\sigma}_{|\sigma|-2}{\longrightarrow}Y^{\sigma}_{|\sigma|-3}{\longrightarrow}\ldots{\longrightarrow}Y^{\sigma}_0$$as follows: $Y^{\sigma}_0=\proj(k[\Xi_i]_{i\in\sigma})\cong\mathbb{P}^{|\sigma|-1}_k$ and $Y_{j+1}^{\sigma}\to Y_{j}^{\sigma}$ is the blowing up of $Y_{j}^{\sigma}$ in the strict transforms (under $Y_{j}^{\sigma}\to Y^{\sigma}_0=\proj(k[\Xi_i]_{i\in\sigma})$) of all $j$-dimensional $k$-rational linear subvarieties of $\proj(k[\Xi_i]_{i\in\sigma})$. We write $Y_{-1}^{\sigma}=\spec(k)$. For $\emptyset\ne\sigma\subset\Upsilon$ let$$H_{\sigma}=\{(a_{ij})_{0\le i,j\le d}\in{\rm GL}\sb {d+1}(k)\quad|\quad\begin{array}{c}a_{ii}=1\,\mbox{for all}\,i\in\sigma^c\\a_{ij}=0\,\mbox{if}\,i\ne j\,\mbox{and}\,\{i,j\}\cap\sigma^{c}\ne\emptyset\end{array}\}.$$Then $H_{\sigma}$-acts on $Y^{\sigma}$ (by forgetting the $a_{ij}$ with $\{i,j\}\cap\sigma^{c}\ne\emptyset$). In several subsequent proofs we will induct on $d$; for that purpose we note that $Y^{\sigma}$ with its action by $H_{\sigma}\cong {\rm GL}\sb {|\sigma|}(k)$ and its fixed ordered set of projective coordinate functions $(\Xi_i)_{i\in\sigma}$ (as ordering on $\sigma$ we take the one induced by its inclusion into the ordered set $\Upsilon=\{0,\ldots,d\}$) is just like the data $$(Y\,\,\mbox{with its}\,\,{\rm GL}\sb {d+1}(k)\,\,\mbox{action},\,\,(\Xi_i)_{0\le i\le d}\,\,\mbox{with the natural ordering on}\,\,\{0,\ldots,d\}),$$ but of dimension $|\sigma|-1$ instead of $d$. For $\tau\subset\sigma$ in ${{{\mathcal{Y}}}}$ let$${{{U}}}_{\tau}^{\sigma}=\{(a_{ij})_{0\le i,j\le d}\in{{{U}}}_{\tau}\quad|\quad a_{ij}=0\,\mbox{if}\,i\ne j\,\mbox{and}\,i\in\sigma^c\}.$$This is a subgroup of $H_{\sigma}$. For $b\in\tau\in{\mathcal Y}$ it is sometimes convenient to write$$U^{\{b\}}_{\tau}:=U^{\Upsilon-\{b\}}_{\tau-\{b\}}=U^{\{b\}^c}_{\tau-\{b\}}.$$For $\sigma\in{{{\mathcal{Y}}}}$ let$${\mathcal N}^{\sigma}=\{(\tau,u)\quad|\quad\emptyset\ne\tau\subsetneq\sigma\,\,\mbox{and}\,\,u\in{{{U}}}^{\sigma}_{\tau}\},$$$${\mathcal V}^{\sigma}=\{V\in{\mathcal V}\quad|\quad V_{\sigma}\ne V\cap V_{{\sigma}}\ne\emptyset\}.$$Then we have a bijection\begin{gather}{\mathcal N}^{\sigma}\coprod{\mathcal N}^{\sigma^c}\cong{\mathcal V}^{\sigma}\label{neigana}\end{gather}where the map ${\mathcal N}^{\sigma}\to{\mathcal V}^{\sigma}$ is given by $(\tau,u)\mapsto u.V_{{\tau}}$, and where the map ${\mathcal N}^{\sigma^c}\to{\mathcal V}^{\sigma}$ is given by $(\tau,u)\mapsto u.V_{{\tau}\cup{\sigma}}$.

\begin{lem}\label{problo} For $\sigma\in{{{\mathcal{Y}}}}$ there is a canonical isomorphism$$V_{\sigma}\cong Y^{\sigma}\times Y^{\sigma^c},$$equivariant for the respective actions of $H_{\sigma}$ and $H_{\sigma^c}$ on both sides.
\end{lem} 

{\sc Proof:} (cf. also \cite{ito}, sect. 4) For $0\le j\le |\sigma|-2$ let $V_{\sigma,j}$ be the closed subscheme of $Y_{j}$ which is the strict transform of $V_{\sigma,0}$ under $Y_{j}\to Y_0$. We can naturally identify $Y^{\sigma^c}$ with $V_{\sigma,d-|\sigma|}$. On the other hand $V_{\sigma,d-|\sigma|+1}=\underline{\proj}(\sym_{{\mathcal O}_{V_{\sigma,d-|\sigma|}}}({\mathcal J}/{\mathcal J}^2))$ according to \cite{hart} II, 8.24, where ${\mathcal J}$ denotes the ideal sheaf of $V_{\sigma,d-|\sigma|}$ in $Y_{d-|\sigma|}$. This is the pull back of the ideal sheaf of $V_{\sigma,0}$ in $Y_0=\proj(k[\Xi_0,\ldots,\Xi_d])$, i.e. the one corresponding to the homogeneous ideal $(\Xi_i)_{i\in\sigma}\subset k[\Xi_0,\ldots,\Xi_d]$. Hence $V_{\sigma,d-|\sigma|+1}=\proj(k[\Xi_i]_{i\in\sigma})\times V_{\sigma,d-|\sigma|}$. Under this isomorphism the successive blowing up of the first factor $\proj(k[\Xi_i]_{i\in\sigma})$ in its $k$-rational linear subvarieties of dimension $\le|\sigma|-3$ corresponds to taking the strict transform of $V_{\sigma,d-|\sigma|+1}$ in $Y$.\hfill$\Box$\\

Let $\sigma\in{\mathcal Y}$. Viewing $Y^{\sigma}$ as the version of $Y$ of dimension $|\sigma|-1$ instead of $d$, definition (\ref{divdef}) provides us with particular divisors on $Y^{\sigma}$. Explicitly we name the divisors$$D(0,{\bf 1},0)^{\sigma}=\sum_{s\in\tau\subsetneq\sigma}\sum_{u\in{{{U}}}_{\tau}^{\sigma}}u.V_{{\tau}}^{\sigma},$$$$D(0,0,{\bf 1})^{\sigma}=\sum_{\emptyset\ne\tau\subset\sigma-\{s\}}\sum_{u\in{{{U}}}_{\tau}^{\sigma}}u.V_{{\tau}}^{\sigma}$$on $Y^{\sigma}$, where $s\in\sigma$ is the minimal element and where the prime divisor $V_{{\tau}}^{\sigma}$ on $Y^{\sigma}$ is the strict transform under $Y^{\sigma}\to Y^{\sigma}_0$ of the common zero set of $\{\Xi_i\}_{i\in\tau}$, and where ${\bf 1}=(1,\ldots,1)\in{\mathbb Z}^{|\sigma|-1}$.

\begin{pro}\label{selfint} Let $\sigma\in{{{\mathcal{Y}}}}$. With the divisor$$E=D(0,{\bf 1},0)^{\sigma}\times Y^{\sigma^c}+Y^{\sigma}\times D(0,0,{\bf 1})^{\sigma^c}$$on $V_{\sigma}=Y^{\sigma}\times Y^{\sigma^c}$ we have the following isomorphism of line bundles on $V_{\sigma}$:$${\mathcal L}_{Y}(-V_{{\sigma}})\otimes_{{\mathcal O}_{Y}}{\mathcal O}_{V_{\sigma}}\cong{\mathcal L}_{V_{\sigma}}(E)$$
\end{pro}

{\sc Proof:} For any $b\in\Upsilon$ the pullback to $Y$ of the divisor $V_{\{b\},0}$ on $Y_0$ is the divisor $\sum_{b\in\tau\in{{{\mathcal{Y}}}}}            \sum_{u\in{{{U}}}_{\tau}^{\{b\}}}u.V_{\tau}$. Let $s$, resp. $t$ be the minimal element of $\sigma$, resp. of $\sigma^c$. The equivalence of divisors $V_{\{s\},0}\sim V_{\{t\},0}$ on $Y_0$ gives rise to the equivalence$$\sum_{s\in\tau\in{{{\mathcal{Y}}}}}         \sum_{u\in{{{U}}}_{\tau}^{\{s\}}}        u.V_{\tau}\quad\sim\quad \sum_{t\in\tau\in{{{\mathcal{Y}}}}}            \sum_{u\in{{{U}}}_{\tau}^{\{t\}}}u.V_{\tau}$$
on $Y$. Thus$$\sum_{s\in\tau\in{{{\mathcal{Y}}}}}       \sum_{u\in{{{U}}}_{\tau}^{\{s\}}\atop u\notin{{{U}}}_{\tau}^{\{t\}}\,\,\mbox{\tiny{if}}\,\,t\in\tau}u.V_{\tau}\quad\sim\quad \sum_{t\in\tau\in{{{\mathcal{Y}}}}}    \sum_{u\in{{{U}}}_{\tau}^{\{t\}}\atop u\notin{{{U}}}_{\tau}^{\{s\}}\,\,\mbox{\tiny{if}}\,\,s\in\tau}    u.V_{\tau}$$
or equivalently$$-V_{\sigma}\quad\sim\quad-(\sum_{t\in\tau\in{{{\mathcal{Y}}}}}    \sum_{u\in{{{U}}}_{\tau}^{\{t\}}\atop u\notin{{{U}}}_{\tau}^{\{s\}}\,\,\mbox{\tiny{if}}\,\,s\in\tau}    u.V_{\tau})+\sum_{s\in\tau\in{{{\mathcal{Y}}}}}     \sum_{u\in{{{U}}}_{\tau}^{\{s\}}\atop {u\notin{{{U}}}_{\tau}^{\{t\}}\,\,\mbox{\tiny{if}}\,\,t\in\tau\atop u\ne1\,\,{\mbox{\tiny{if}}}\,\,\tau=\sigma}}u.V_{\tau}.$$ Now we are interested only in the summands which belong to the set ${\mathcal V}^{\sigma}$ which we determined in (\ref{neigana}). For example, all the summands in the bracketed term on the right hand side do not belong to ${\mathcal V}^{\sigma}$ (the condition $t\in\tau$ excludes contributions from ${\mathcal N}^{\sigma}$, the condition $u\notin{{{U}}}_{\tau}^{\{s\}}\,\mbox{if}\,s\in\tau$ and the fact ${{{U}}}_{\tau-\sigma}^{\sigma^c}\subset{{{U}}}_{\tau}^{\{s\}}$ for $\sigma\subset\tau$ exclude contributions from ${\mathcal N}^{\sigma^c}$). We get$${\mathcal L}_{Y}(-V_{\sigma})\otimes_{{\mathcal O}_{Y}}{\mathcal O}_{V_{\sigma}}\cong{\mathcal L}_{V_{\sigma}}(E)$$with $$E=V_{\sigma}\cap(\sum_{{s\in\tau\subsetneq\sigma}}\sum_{u\in{{{U}}}_{\tau}^{\{s\}}\cap{{{U}}}_{\tau}^{\sigma}}u.V_{\tau}+\sum_{{\sigma\subsetneq\tau\subsetneq\Upsilon}}\sum_{u\in{{{U}}}_{\tau-\sigma}^{\sigma^c}\atop u\notin{{{U}}}_{\tau}^{\{t\}}\,\,\mbox{\tiny{if}}\,\,t\in\tau}u.V_{\tau}).$$ Now since $s$ is minimal in $\sigma$ we see that ${{{U}}}^{\sigma}_{\tau}\cap {{{U}}}^{\{s\}}_{\tau}={{{U}}}^{\sigma}_{\tau}$ for all $\tau\subset\sigma$. On the other hand, since $t$ is minimal in $\sigma^c$ we have ${{{U}}}_{\tau}^{\{t\}}={{{U}}}_{\tau}$ for all $\sigma\subset\tau$ with $t\in\tau$. We obtain$$E=V_{\sigma}\cap(\sum_{{s\in\tau\subsetneq\sigma}}\sum_{u\in{{{U}}}_{\tau}^{\sigma}}u.V_{\tau}+\sum_{{\sigma\subsetneq\tau\subset\Upsilon-\{t\}}}\sum_{{u\in{{{U}}}_{\tau-\sigma}^{\sigma^c}}}u.V_{\tau}).$$That this is the divisor $E$ as stated follows from the construction of the isomorphism $V_{\sigma}=Y^{\sigma}\times Y^{\sigma^c}$.\hfill$\Box$\\

\begin{pro}\label{nullsta} Assume $-d\le n_1\le\ldots\le n_d$ and $0\le m_1\le \ldots\le m_d$, furthermore $n_i-n_1\le i-1$ and $m_i-m_1\le i-1$ for all $1\le i\le d$.\\(a) We have $H^t(Y,{\mathcal L}_Y(D(0,n,m)))=0$ for all $t\in \mathbb{Z}_{>0}$.\\(b) If $n_d\le -1$ and $m_1=0$ we have $H^t(Y,{\mathcal L}_Y(D(0,n,m)))=0$ for all $t\in \mathbb{Z}$.\\
\end{pro}

{\sc Proof:} (b) Outer induction on $d$, inner induction on $$s(n,m)=\sum_{i=1}^d(m_i+n_i-n_1).$$Our assumptions imply $s(n,m)\ge0$. If $s(n,m)=0$ or if $d=1$ we have $n_i=n_1$ and $m_i=0$ for all $1\le i\le d$. Then $D(0,n,m)$ is the pullback of ${\mathcal O}_{\mathbb{P}^{d}}(n_1)$ under the successive blowing up $Y\to\mathbb{P}^{d}$ and the claim follows from $H^t(\mathbb{P}^{d},{\mathcal O}_{\mathbb{P}^{d}}(n_1))=0$ for all $t\in \mathbb{Z}$. Now let $s(n,m)>0$. First suppose that there exists a $2\le i\le d$ with $n_i\ne n_1$. Then let $i_0$ be minimal with this property. Let $n'_{i_0}=n_{i_0}-1$ and $n'_i=n_i$ for all $i\ne i_0$. Then also $n'=(n'_1,\ldots,n'_d)$ satisfies our hypothesis, and $s(n',m)<s(n,m)$. We have an exact sequence$$0\longrightarrow{\mathcal L}_Y(D(0,n',m))\longrightarrow{\mathcal L}_Y(D(0,n,m))\longrightarrow{\mathcal C}\longrightarrow0$$and in view of the induction hypothesis it suffices to show $H^t(Y,{\mathcal C})=0$ for all $t\in \mathbb{Z}$. We may view ${\mathcal C}$ as living on the closed subscheme $\coprod_{(u,\sigma)}u.V_{\sigma}$ with $\sigma$ running through all elements of ${{{\mathcal{Y}}}}$ with $0\in\sigma$ and $|\sigma|=i_0$, and with $u\in{{{U}}}_{\sigma}$. We deal with every such $(u,\sigma)$ separately. By equivariance we may assume $u=1$. The restriction ${\mathcal C}|_{V_{\sigma}}$ of ${\mathcal C}$ to $V_{\sigma}=Y^{\sigma}\times Y^{\sigma^c}$ is isomorphic to$${\mathcal L}_{Y^{\sigma}}(D(0,\widetilde{n},\widetilde{m}))\otimes_{{\mathcal O}_{Y^{\sigma}}}{\mathcal O}_{V_{\sigma}}\otimes_{{\mathcal O}_{Y^{\sigma^c}}}{\mathcal E}$$with a line bundle ${\mathcal E}$ on $Y^{\sigma^c}$ and with $\widetilde{n}=(\widetilde{n}_1,\ldots,\widetilde{n}_{i_0-1})$ and $\widetilde{m}=(\widetilde{m}_1,\ldots,\widetilde{m}_{i_0-1})$ defined by $\widetilde{n}_i=n_{i}-n_{i_0}$ and $\widetilde{m}_i=m_i$ (for $1\le i\le i_{0}-1$). This follows from Proposition \ref{selfint}. By induction hypothesis we have $H^t(Y^{\sigma},{\mathcal L}_{Y^{\sigma}}(D(0,\widetilde{n},\widetilde{m})))=0$ for all $t\in \mathbb{Z}$. Thus $H^t(Y,{\mathcal C}|_{V_{\sigma}})=0$ follows from the K\"unneth formula. 

If there is no $2\le i\le d$ with $n_i\ne n_1$ then $s(n,m)>0$ implies that there is a $2\le i\le d$ with $m_i\ne0$. Let $i_0$ be minimal with this property. Let $m'_{i_0}=m_{i_0}-1$ and $m'_i=m_i$ for all $i\ne i_0$. Then also $m'=(m'_1,\ldots,m'_d)$ satisfies our hypothesis, and $s(n,m')<s(n,m)$. We have an exact sequence$$0\longrightarrow{\mathcal L}_Y(D(0,n,m'))\longrightarrow{\mathcal L}_Y(D(0,n,m))\longrightarrow{\mathcal C}\longrightarrow0$$and in view of the induction hypothesis it suffices to show $H^t(Y,{\mathcal C})=0$ for all $t\in \mathbb{Z}$. We may view ${\mathcal C}$ as living on the closed subscheme $\coprod_{(u,\sigma)}u.V_{\sigma}$ with $\sigma$ running through all elements of ${{{\mathcal{Y}}}}$ with $0\notin\sigma$ and $|\sigma|=i_0$, and with $u\in{{{U}}}_{\sigma}$. We deal with every such $(u,\sigma)$ separately. By equivariance we may assume $u=1$. The restriction ${\mathcal C}|_{V_{\sigma}}$ of ${\mathcal C}$ to $V_{\sigma}=Y^{\sigma}\times Y^{\sigma^c}$ is isomorphic to$${\mathcal L}_{Y^{\sigma}}(D(0,\widetilde{n},\widetilde{m}))\otimes_{{\mathcal O}_{Y^{\sigma}}}{\mathcal O}_{V_{\sigma}}\otimes_{{\mathcal O}_{Y^{\sigma^c}}}{\mathcal E}$$with a line bundle ${\mathcal E}$ on $Y^{\sigma^c}$ and with $\widetilde{n}=(\widetilde{n}_1,\ldots,\widetilde{n}_{i_0-1})$ and $\widetilde{m}=(\widetilde{m}_1,\ldots,\widetilde{m}_{i_0-1})$ defined by $\widetilde{n}_i=m_{i}-m_{i_0}$ and $\widetilde{m}_i=m_i$ (for $1\le i\le i_{0}-1$). This follows from Proposition \ref{selfint}. By induction hypothesis we have $H^t(Y^{\sigma},{\mathcal L}_{Y^{\sigma}}(D(0,\widetilde{n},\widetilde{m})))=0$ for all $t\in \mathbb{Z}$. Thus $H^t(Y,{\mathcal C}|_{V_{\sigma}})=0$ follows from the K\"unneth formula.

(a) Outer induction on $d$. If $d=1$ our statement is the well known fact $H^t(\mathbb{P}^{1},{\mathcal O}_{\mathbb{P}^{1}}(k))=0$ for all $k\in\mathbb{Z}_{\ge-1}$, all $t\in \mathbb{Z}_{>0}$. Inner induction on$$r(n,m)=d^2+\sum_{i=1}^d(m_i+n_i)$$which is $\ge0$ as follows from our assumptions. The case $r(n,m)=0$ corresponds to $n_i=-d$ and $m_i=0$ for all $i$, hence was settled in (b). Now suppose $r(n,m)>0$. Then at least one of the following cases occurs:\\(i) there is a $1\le i_0\le d$ such that if we set $m'_i=m_i$ for $i\ne i_0$ and $m'_{i_0}=m_{i_0}-1$, then also $m'=(m'_1,\ldots,m'_d)$ satisfies our hypothesis.\\(ii)  there is a $1\le i_0\le d$ such that if we set $n'_i=n_i$ for $i\ne i_0$ and $n'_{i_0}=n_{i_0}-1$, then also $n'=(n'_1,\ldots,n'_d)$ satisfies our hypothesis.\\Fix one of these two cases which holds true and let $i_0$ be the minimal element satisfying its condition. In case (i) let $n'=n$ and in case (ii) let $m'=m$. Since by induction hypothesis we have $H^t(Y,{\mathcal L}_Y(D(0,n',m')))=0$ for all $t\in \mathbb{Z}_{>0}$, the exact sequence$$0\longrightarrow{\mathcal L}_Y(D(0,n',m'))\longrightarrow{\mathcal L}_Y(D(0,n,m))\longrightarrow{\mathcal C}\longrightarrow0$$shows that it suffices to show $H^t(Y,{\mathcal C})=0$ for all $t\in \mathbb{Z}_{>0}$. We may view ${\mathcal C}$ as living on the closed subscheme $\coprod_{(u,\sigma)}u.V_{\sigma}$ with $u\in{{{U}}}_{\sigma}$ and with $\sigma$ running through all elements of ${{{\mathcal{Y}}}}$ with $|\sigma|=i_0$ and in addition with $0\in\sigma$ in case (ii), resp. with $0\notin\sigma$ in case (i). We deal with every such $(u,\sigma)$ separately. By equivariance we may assume $u=1$. The restriction ${\mathcal C}|_{V_{\sigma}}$ of ${\mathcal C}$ to $V_{\sigma}=Y^{\sigma}\times Y^{\sigma^c}$ is isomorphic (use Proposition \ref{selfint}) to$${\mathcal L}_{Y^{\sigma}}(D(0,\widetilde{n},\widetilde{m}))\otimes_{{\mathcal O}_{Y^{\sigma}}}{\mathcal O}_{V_{\sigma}}\otimes_{{\mathcal O}_{Y^{\sigma^c}}}{\mathcal L}_{Y^{\sigma^c}}(D(0,\widehat{n},\widehat{m}))$$with $\widetilde{m}=(\widetilde{m}_1,\ldots,\widetilde{m}_{i_0-1})$, $\widetilde{n}=(\widetilde{n}_1,\ldots,\widetilde{n}_{i_0-1})$, $\widehat{m}=(\widehat{m}_1,\ldots,\widehat{m}_{d-i_0})$ and $\widehat{n}=(\widehat{n}_1,\ldots,\widehat{n}_{d-i_0})$ defined as follows: $\widehat{m}_i=m_{i+i_0}-m_{i_0}$ in case (i), and $\widehat{m}_i=n_{i+i_0}-n_{i_0}$ in case (ii). Moreover $\widetilde{n}_i=m_i-m_{i_0}$ in case (i), and $\widetilde{n}_i=n_i-n_{i_0}$ in case (ii). Finally $\widetilde{m}_i=m_i$ and $\widehat{n}_i=n_{i_0+i}$ in all cases. By induction hypothesis $H^t(Y^{\sigma},{\mathcal L}_{Y^{\sigma}}(D(0,\widetilde{n},\widetilde{m})))=0$ and $H^t(Y^{\sigma^c},{\mathcal L}_{Y^{\sigma^c}}(D(0,\widehat{n},\widehat{m})))=0$ for all $t\in \mathbb{Z}_{>0}$ so we conclude by the K\"unneth formula.\hfill$\Box$\\

\begin{satz}\label{gallgvan} Assume $-d\le n_1\le\ldots\le n_d$ and $0\le m_1\le \ldots\le m_d$, furthermore $n_{i+1}-n_i\le 1$ and $m_{i+1}-m_i\le 1$ for all $1\le i\le d-1$. Let$$e=\min\{0\le i\le d\quad|\quad \overline{a}_t\le0 \mbox{ for all } t>i\}.$$ Then $H^t(Y,{\mathcal L}_Y(D(\overline{a},n,m)))=0$ for all $t>e$.
\end{satz}  

{\sc Proof:} Let us introduce some more notation. For an element $\sigma\in{{{\mathcal{Y}}}}$ we define $$\overline{a}|_{\sigma}=(\overline{a}|_{\sigma,1},\ldots,\overline{a}|_{\sigma,|\sigma|-1})\in\mathbb{Z}^{|\sigma|-1}$$as follows. Let $$\iota_{\sigma}:\{0,\ldots,|\sigma|-1\}\to\sigma$$ be the order preserving bijection. Then$$\overline{a}|_{\sigma,i}=\overline{a}_{\iota_{\sigma}(i)}\quad\quad(1\le i\le|\sigma|-1).$$Thus, $\overline{a}|_{\sigma}$ enumerates those components of the $(d+1)$-tupel $(\overline{a}_0,\ldots,\overline{a}_d)$ which have index in $\sigma$, but omits the first one of these. For an $s$-tuple $\overline{b}=(\overline{b}_1,\ldots,\overline{b}_s)\in\mathbb{Z}^s$ (some $s\in\mathbb{N}$) and an element $i\in\{1,\ldots,s\}$ we define $$\overline{b}^{[i]}=(\overline{b}^{[i]}_{1},\ldots,\overline{b}^{[i]}_{s})=\overline{b}+(0,\ldots,0,1,0,\ldots,0)\in\mathbb{Z}^{s}$$by setting $\overline{b}^{[i]}_{j}=\overline{b}_{j}$ for $j\ne i$, and $\overline{b}^{[i]}_{i}=\overline{b}_{i}+1$. Moreover we let $\overline{b}^{[0]}=\overline{b}$.\\

Now we begin. Outer induction on $d$. If $d=1$ our statement is the well known fact $H^t(\mathbb{P}^{1},{\mathcal O}_{\mathbb{P}^{1}}(k))=0$ for all $k\in\mathbb{Z}_{\ge-1}$, all $t\in \mathbb{Z}_{>0}$. For $d>1$ we proceed in two steps.\\

{\it  First Step: The case $\overline{a}_i\ge 0$ for all $1\le i\le d$.}\\ Induction on $-\overline{a}_0=\sum_{i=1}^d\overline{a}_i$. The case $-\overline{a}_0=0$ is settled in Proposition \ref{nullsta}. (If $e=0$ --- the case relevant for our later Theorems \ref{wellvan} and \ref{globnull} --- this is the end of the present first step.) For the induction step suppose we are given a $d$-tupel $\overline{a}=(\overline{a}_1,\ldots,\overline{a}_d)$ with $\overline{a}_t\le 0$ for all $t>e$ and an integer $1\le j\le e$ with $\overline{a}_j\ge0$. Suppose we know $H^t(Y,{\mathcal L}_Y(D(\overline{a},n,m)))=0$ for all $t>e$ (induction hypothesis). Defining $\overline{a}'=(\overline{a}'_1,\ldots,\overline{a}'_d)$ by $\overline{a}'_i=\overline{a}_i$ for $i\ne j$, and $\overline{a}'_j=\overline{a}_j+1$, we need to show $H^t(Y,{\mathcal L}_Y(D(\overline{a}',n,m)))=0$ for all $t>e$. For $0\le k\le d$ let \begin{gather}D_k(\overline{a},n,m)=D(\overline{a},n,m)-\sum_{j\in\tau\in{\mathcal Y}\atop |\tau|\le k}\sum_{u\in{{{U}}}_{\tau}-{{{U}}}_{\tau}^{\{j\}}}u.V_{\tau}\label{deka}.\end{gather}Thus $D(\overline{a},n,m)=D_0(\overline{a},n,m)$. On the other hand $D_d(\overline{a},n,m)\sim D(\overline{a}',n,m)$. Indeed, $$D(\overline{a}',n,m)=D(\overline{a},n,m)+\sum_{0\in\tau\in{\mathcal Y}}\sum_{u\in U_{\tau}}u.V_{\tau}-\sum_{j\in\tau\in{\mathcal Y}}\sum_{u\in U_{\tau}}u.V_{\tau}$$and we have $$\sum_{0\in\tau\in{\mathcal Y}}\sum_{u\in U_{\tau}}u.V_{\tau}\sim \sum_{j\in\tau\in{\mathcal Y}}\sum_{u\in U_{\tau}^{\{j\}}}u.V_{\tau}$$on $Y$ because $\sum_{0\in\tau\in{\mathcal Y}}\sum_{u\in U_{\tau}}u.V_{\tau}-\sum_{j\in\tau\in{\mathcal Y}}\sum_{u\in U_{\tau}^{\{j\}}}u.V_{\tau}$ is the (principal) divisor of the rational function $\Xi_0/\Xi_j$ on $Y$ (note that $U_{\tau}=U_{\tau}^{\{0\}}$ if $0\in\tau$). Therefore, in view of our induction hypothesis and of the exact sequence$$0\longrightarrow{\mathcal L}_Y(D_d(\overline{a},n,m))\longrightarrow {\mathcal L}_Y(D_0(\overline{a},n,m))\longrightarrow\frac{{\mathcal L}_Y(D_0(\overline{a},n,m))}{{\mathcal L}_Y(D_d(\overline{a},n,m))}\longrightarrow0$$we only need to show $$H^t(Y,\frac{{\mathcal L}_Y(D_0(\overline{a},n,m))}{{\mathcal L}_Y(D_d(\overline{a},n,m))})=0$$ for all $t\ge e$. By the obvious induction we reduce to proving $$H^t(Y,\frac{{\mathcal L}_Y(D_{k-1}(\overline{a},n,m))}{{\mathcal L}_Y(D_{k}(\overline{a},n,m))})=0$$ for all $t\ge e$, all $1\le k\le d$. We have$$\frac{{\mathcal L}_Y(D_{k-1}(\overline{a},n,m))}{{\mathcal L}_Y(D_{k}(\overline{a},n,m))}=\bigoplus_{j\in\tau\in{\mathcal Y}\atop |\tau|=k}\bigoplus_{u\in{{{U}}}_{\tau}-{{{U}}}_{\tau}^{\{j\}}}{\mathcal L}_Y(D_{k-1}(\overline{a},n,m))\otimes_{{\mathcal O}_Y}{\mathcal O}_{u.V_{\tau}}$$so our task is to prove $$H^t(Y,{\mathcal L}_Y({D}_{k-1}(\overline{a},n,m))\otimes_{{\mathcal O}_Y}{\mathcal O}_{u.V_{\tau}})=0$$for any $t\ge e$, any $\tau$ with $j\in\tau$ and $|\tau|=k$,
and any $u\in{{{U}}}_{\tau}-{{{U}}}_{\tau}^{\{j\}}$. Setting\begin{gather}\widehat{D}_{k-1}(\overline{a},n,m)=D(\overline{a},n,m)-\sum_{j\in\sigma\in{\mathcal Y}\atop |\sigma|\le k-1}\sum_{u\in{{{U}}}_{\sigma}}u.V_{\sigma}\label{dehut}\end{gather}we claim $${\mathcal L}_Y({D}_{k-1}(\overline{a},n,m))\otimes_{{\mathcal O}_Y}{\mathcal O}_{u.V_{\tau}}={\mathcal L}_Y(\widehat{D}_{k-1}(\overline{a},n,m))\otimes_{{\mathcal O}_Y}{\mathcal O}_{u.V_{\tau}}.$$Indeed, to show this we need to show that for all $\tau'$ with $j\in\tau'\subsetneq\tau$, all $u'\in {{{U}}}_{\tau'}^{\{j\}}$ we have $u'.V_{\tau'}\cap u.V_{\tau}=\emptyset$. Assume that this is false. Then for the images $f(V_{\tau'})=\mathbb{V}(\Xi_i)_{i\in\tau'}$ and $f(V_{\tau})=\mathbb{V}(\Xi_i)_{i\in\tau}$ under $f:Y\to Y_0=\proj(k[\Xi_0,\ldots,\Xi_d])$ (we use the symbol ${\mathbb{V}}(.)$ to denote the vanishing locus of a set of projective coordinate functions) we have $$u.\mathbb{V}(\Xi_i)_{i\in\tau}\subset u.\mathbb{V}(\Xi_i)_{i\in\tau'}$$ or equivalently$$(u')^{-1}u.\mathbb{V}(\Xi_i)_{i\in\tau}\subset\mathbb{V}(\Xi_i)_{i\in\tau'}.$$Write $u^{-1}=(c_{qr})_{qr}$ and $u'=(d_{qr})_{qr}$ and $u^{-1}u'=(e_{qr})_{qr}$. Because of $u^{-1}\in{{{U}}}_{\tau}-{{{U}}}_{\tau}^{\{j\}}$ there is a $s<j$, $s\notin\tau$, such that $c_{sj}\ne0$. On the other hand $d_{qj}=0$ for all $q\ne j$. Thus $e_{sj}=c_{sj}\ne0$. But then if $P\in\mathbb{V}(\Xi_i)_{i\in\tau}$ denotes the point with homogeneous coordinates $=0$ for all indices $\ne s$ we have $((u')^{-1}u)P\notin\mathbb{V}(\Xi_j)$, in particular $((u')^{-1}u)P\notin\mathbb{V}(\Xi_i)_{i\in\tau'}$. Thus the assumption was false and the claim is established. Now notice that $\widehat{D}_{k-1}(\overline{a},n,m)$ is invariant under ${{{U}}}$, implying that under the isomorphism $u:V_{\tau}\to u.V_{\tau}$ we have $${\mathcal L}_Y(\widehat{D}_{k-1}(\overline{a},n,m))\otimes_{{\mathcal O}_Y}{\mathcal O}_{u.V_{\tau}}\cong{\mathcal L}_Y(\widehat{D}_{k-1}(\overline{a},n,m))\otimes_{{\mathcal O}_Y}{\mathcal O}_{V_{\tau}}.$$In this way we have transformed our task into that of proving$$H^t(Y,{\mathcal L}_Y(\widehat{D}_{k-1}(\overline{a},n,m))\otimes_{{\mathcal O}_Y}{\mathcal O}_{V_{\tau}})=0$$for any $t\ge e$ and any $\tau$ with $j\in\tau$ and $|\tau|=k$ and $U_{\tau}-U_{\tau}^{\{j\}}\ne\emptyset$. We use Proposition \ref{selfint} to show that on $V_{\tau}=Y^{\tau}\times Y^{\tau^c}$ we have \begin{gather}{\mathcal L}_Y(\widehat{D}_{k-1}(\overline{a},n,m))\otimes_{{\mathcal O}_Y}{\mathcal O}_{V_{\tau}}\notag\cong\\{\mathcal L}_{Y^{\tau}}({D}((\overline{a}|_{\tau})^{[\iota_{\tau}^{-1}(j)]},n',m'))\otimes_{{\mathcal O}_{Y^{\tau}}}{\mathcal O}_{V_{\tau}}\otimes_{{\mathcal O}_{Y^{\tau^c}}}{\mathcal L}_{Y^{\tau^c}}({D}(\overline{a}|_{\tau^c},n'',m''))\label{tepro}\end{gather}with ${n}'=({n}'_1,\ldots,{n}'_{k-1})$, ${m}'=({m}'_1,\ldots,{m}'_{k-1})$, ${n}''=({n}''_1,\ldots,{n}''_{d-k})$ and ${m}''=({m}''_1,\ldots,{m}''_{d-k})$ defined as follows: $$m_i'=m_i\,\,\mbox{for all}\,\,1\le i\le k-1\quad\quad n_i''=n_{k+i}\,\,\mbox{for all}\,\,1\le i\le d-k$$$$n_i'=\left\{\begin{array}{l@{\quad:\quad}l}m_{i}-m_k-1&0\notin\tau\\n_{i}-n_k-1&0\in\tau\end{array}\right.\quad\quad m_i''=\left\{\begin{array}{l@{\quad:\quad}l}m_{k+i}-m_k&0\notin\tau\\n_{k+i}-n_k&0\in\tau\end{array}\right.$$By induction hypothesis we have $$H^t(Y^{\tau},{\mathcal L}_{Y^{\tau}}({D}((\overline{a}|_{\tau})^{[\iota_{\tau}^{-1}(j)]},n',m')))=0$$for all $t>e'$, and$$H^t(Y^{\tau^c},{\mathcal L}_{Y^{\tau^c}}({D}(\overline{a}|_{\tau^c},n'',m'')))=0$$ for all $t>e''$. Here $e'$ and $e''$ are defined as follows. For a non empty subset $\sigma$ of $\{0,\ldots,d\}$ denote by $\widehat{\sigma}$ the subset of $\sigma$ obtained by deleting its smallest element, and all its elements $>e$. Then $e'=|\widehat{{\tau}}|$ and $e''=|\widehat{{\tau^c}}|$. Now observe $e'+e''<e$. Indeed, otherwise we would either have $\{0,\ldots,e\}\subset\tau^c$ --- contradicting $j\in\{0,\ldots,e\}\cap \tau$ --- or we would have $\{0,\ldots,e\}\subset\tau$ --- contradicting $j\in\{0,\ldots,e\}$ and $U_{\tau}-U_{\tau}^{\{j\}}\ne\emptyset$. We conclude this step by the K\"unneth formula.\\

{\it  Second Step: The general case.}\\Induction on $-\sum_{\overline{a}_j<0}\overline{a}_j$. The case where this term is zero was settled in the first step. Now let $-\sum_{\overline{a}_j<0}\overline{a}_j>0$. Choose and fix a $j$ with $\overline{a}_j<0$. Define $\overline{a}'=(\overline{a}'_1,\ldots,\overline{a}'_d)$ by $\overline{a}'_i=\overline{a}_i$ for $i\ne j$, and $\overline{a}'_j=\overline{a}_j+1$. For $0\le k\le d$ define $D_k(\overline{a},n,m)$ as in the first step through formula (\ref{deka}). Then we have again $D_d(\overline{a},n,m)\sim D(\overline{a}',n,m)$ and similarly to the first step we use the induction hypothesis to reduce to proving$$H^t(Y,\frac{{\mathcal L}_Y(D_{k-1}(\overline{a},n,m))}{{\mathcal L}_Y(D_{k}(\overline{a},n,m))})=0$$ for all $t>e$, all $1\le k\le d$ (note that now we write $t>e$ rather than $t\ge e$). Similarly to the first step we then reduce to proving$$H^t(Y,{\mathcal L}_Y({D}_{k-1}(\overline{a},n,m))\otimes_{{\mathcal O}_Y}{\mathcal O}_{u.V_{\tau}})=0$$for any $t>e$, any $\tau$ with $j\in\tau$ and $|\tau|=k$,
and any $u\in{{{U}}}_{\tau}-{{{U}}}_{\tau}^{\{j\}}$. Again we define $\widehat{D}_{k-1}(\overline{a},n,m)$ through formula (\ref{dehut}) and similarly as in ther first case we reduce to proving$$H^t(Y,{\mathcal L}_Y(\widehat{D}_{k-1}(\overline{a},n,m))\otimes_{{\mathcal O}_Y}{\mathcal O}_{V_{\tau}})=0$$for any $t>e$ and any $\tau$ with $j\in\tau$ and $|\tau|=k$ and $U_{\tau}-U_{\tau}^{\{j\}}\ne\emptyset$. Again we have the formula (\ref{tepro}) and by induction hypothesis we have $$H^t(Y^{\tau},{\mathcal L}_{Y^{\tau}}({D}((\overline{a}|_{\tau})^{[\iota_{\tau}^{-1}(j)]},n',m')))=0$$for all $t>e'$, and$$H^t(Y^{\tau^c},{\mathcal L}_{Y^{\tau^c}}({D}(\overline{a}|_{\tau^c},n'',m'')))=0$$ for all $t>e''$, where $e'$ and $e''$ are defined as in the first step. Now we clearly have $e'+e''\le e$ (but not necessarily $e'+e''< e$ as in the first case) and we conclude by the K\"unneth formula.\hfill$\Box$\\

For an element of ${\mathcal V}$ let us denote by a subscript $0$ its image in $Y_0\cong {\mathbb P}^d$, i.e. the $k$-linear subspace of $Y_0$ whose strict transform is the given element of ${\mathcal V}$. Let $V, V'\in {\mathcal V}$ and suppose that $V_0\cup V'_0$ is contained in a {\it proper} $k$-linear subspace of $Y_0$; then denote by $[V,V']$ the element of ${\mathcal V}$ which is the strict transform of the minimal such subspace of $Y_0$. If $V_0\cup V'_0$ is not contained in a proper $k$-linear subspace of $Y_0$ than $[V,V']$ is undefined. We say that a subset $S$ of ${\mathcal V}$ is {\it stable} if for any two $V, V'\in S$ the element $[V,V']$ is defined and lies in $S$. For example, the empty set and all one-element subsets of ${\mathcal V}$ are stable.

For any subset $S$ of ${\mathcal V}$ we define the divisor $$D_S=\sum_{V\in S}V$$on $Y$.

\begin{satz}\label{wellvan} Suppose that $S$ is stable and that $\overline{a}=(\overline{a}_1,\ldots,\overline{a}_d)\in\mathbb{Z}^d$ satisfies $\overline{a}_i\le0$ for all $1\le i\le d$. Then $$H^t(Y,{\mathcal L}_Y(D(\overline{a},0,0)-D_S))=0$$ for all $t\in \mathbb{Z}_{>0}$.
\end{satz}

{\sc Proof:} Outer induction on $d$. If $d=1$ our statement is the well known fact $H^t(\mathbb{P}^{1},{\mathcal O}_{\mathbb{P}^{1}}(k))=0$ for all $k\in\mathbb{Z}_{\ge-1}$, all $t\in \mathbb{Z}_{>0}$. Let now $d>1$. For an element $V\in{\mathcal V}$ let$$T(V)=\{V'\in {\mathcal V}|\quad V'_0\subset V_0\}.$$If $S=\emptyset$ we cite Theorem \ref{gallgvan}. Otherwise there is a $V\in S$ with $S\subset T(V)$. It follows that there is a $W\in{\mathcal V}$ (not necessarily $W\in S$) with $\dim(W)=d-1$ and $S\subset T(W)$. For $1\le i\le d$ set $$Q_i=S\cup\{V\in T(W)|\,\,\dim(V_0)\ge i\}.$$By induction on $i$ we show$$H^t(Y,{\mathcal L}_Y(D(\overline{a},0,0)-D_{Q_i}))=0$$ for all $t\in \mathbb{Z}_{>0}$; for $i=d$ we then get our claim. For $i=0$ we have $Q_0=T(W)$ so we see that the divisor $D(\overline{a},0,0)-D_{Q_0}$ is linearly equivalent to $D(\overline{a},{{\bf{-1}}},0)$ with ${{\bf{-1}}}=(-1,\ldots,-1)$. Thus Theorem \ref{gallgvan} settles this case. For $i>0$ we have an exact sequence$$0\longrightarrow{\mathcal L}_Y(D(\overline{a},0,0)-D_{Q_{i-1}})\longrightarrow {\mathcal L}_Y(D(\overline{a},0,0)-D_{Q_i})\longrightarrow{\mathcal C}\longrightarrow0$$with$${\mathcal C}=\bigoplus_{V\in T(W)-S\atop \dim(V_0)=i-1} {\mathcal L}_Y(D(\overline{a},0,0)-D_{Q_i})\otimes_{{\mathcal O}_Y}{\mathcal O}_{V}.$$In view of the induction hypothesis it remains to show $$H^t(Y,{\mathcal L}_Y(D(\overline{a},0,0)-D_{Q_i})\otimes_{{\mathcal O}_Y}{\mathcal O}_{V})=0$$ for all $t\in \mathbb{Z}_{>0}$, all $V\in T(W)-S$ with $\dim(V_0)=i-1$. By equivariance we may assume $V=V_{\tau}$ with $|\tau|=d+1-i$. Recall the bijection (\ref{neigana}) which we now view as an identification as follows: $${\mathcal N}^{\tau^c}=\{V\in{\mathcal V}\,\,|\,\,V_0\subsetneq V_{\tau,0}\},$$ $${\mathcal N}^{\tau}=\{V\in{\mathcal V}\,\,|\,\,V_{\tau,0}\subsetneq V_0\}.$$ On $V_{\tau}=Y^{\tau}\times Y^{\tau^c}$ we find (using Proposition \ref{selfint}) $${\mathcal L}_Y(D(\overline{a},0,0)-D_{Q_i})\otimes_{{\mathcal O}_Y}{\mathcal O}_{V_{\tau}}\cong$$$${\mathcal L}_{Y^{\tau}}(D(\overline{a}|_{\tau},0,0)-D_{T(W)\cap {\mathcal N}^{\tau}})\otimes_{{\mathcal O}_{Y^{\tau}}}{\mathcal O}_{V_{\tau}}\otimes_{{\mathcal O}_{Y^{\tau^c}}}{\mathcal L}_{Y^{\tau^c}}(D(\overline{a}|_{\tau^c},0,0)-D_{S\cap {\mathcal N}^{\tau^c}})$$(with $\overline{a}|_{\tau}$ and $\overline{a}|_{\tau^c}$ as defined in the proof of Theorem \ref{gallgvan}). Since $S$ is stable (with respect to ${\mathcal V}\cong{\mathcal N}$) also $S\cap {\mathcal N}^{\tau^c}$ is stable with respect to ${\mathcal N}^{\tau^c}$. On the other hand $T(W)\cap {\mathcal N}^{\tau}$ is stable with respect to ${\mathcal N}^{\tau}$. Therefore our induction hypothesis says $H^t(Y^{\tau},{\mathcal L}_{Y^{\tau}}(D(\overline{a}|_{\tau},0,0)-D_{T(W)\cap {\mathcal N}^{\tau}}))=0$ and $H^t(Y^{\tau^c},{\mathcal L}_{Y^{\tau^c}}(D(\overline{a}|_{\tau^c},0,0)-D_{S\cap {\mathcal N}^{\tau^c}}))=0$ for all $t\in \mathbb{Z}_{>0}$. We conclude by the K\"unneth formula.\hfill$\Box$\\

For subsets $S$ of ${\mathcal V}$ let $W_S=\bigcap_{V\in S}V$; in particular let $W_{\emptyset}=Y$.

\begin{kor}\label{kritvan} Suppose $\overline{a}$ is as in Theorem \ref{wellvan}.\\(a) Let $S,S'$ be subsets of ${\mathcal V}$ (possibly empty) with $W_{S\cup S'}\ne\emptyset$ and $S\cap S'=\emptyset$. Then $$H^t(Y,{\mathcal L}_Y(D(\overline{a},0,0)-D_{S'})\otimes_{{\mathcal O}_Y}{\mathcal O}_{W_S})=0$$ for all $t\in \mathbb{Z}_{>0}$.\\(b) Let $S$ be a subset of ${\mathcal V}$ and let$$N=\{V\in\mathcal{V}|\,\,V'_0\subset V_0 \,\,\mbox{ for all }\,\,V'\in S\}.$$Let $M_0$ be a non empty and stable subset of $N$. Then the sequence$$H^0(Y,{\mathcal L}_Y(D(\overline{a},0,0))\otimes_{{\mathcal O}_Y}{\mathcal O}_{W_S})\longrightarrow\prod_{V\in M_0}H^0(Y,{\mathcal L}_Y(D(\overline{a},0,0))\otimes_{{\mathcal O}_Y}{\mathcal O}_{W_{S\cup\{V\}}})$$$$\longrightarrow\prod_{V\ne V'\in M_0\atop W_{S\cup\{V,V'\}}\ne\emptyset}H^0(Y,{\mathcal L}_Y(D(\overline{a},0,0))\otimes_{{\mathcal O}_Y}{\mathcal O}_{W_{S\cup\{V,V'\}}})$$is exact.
\end{kor}

{\sc Proof:} (a) The condition $W_S\cap W_{S'}\ne \emptyset$ implies that all subsets of $S\cup S'$ are stable. Therefore we can use Theorem \ref{wellvan} for an induction on $|S|$.\\(b) The kernel of the second arrow is $H^0(Y,{\mathcal L}_Y(D(\overline{a},0,0))\otimes_{{\mathcal O}_Y}{\mathcal O}_{W})$ with $W=\bigcup_{V\in M_0}W_{S\cup \{V\}}$ (union inside $Y$). Looking at the exact sequence$$0\longrightarrow{\mathcal L}_Y(D(\overline{a},0,0)-D_{M_0})\otimes_{{\mathcal O}_Y}{\mathcal O}_{W_S}\longrightarrow{\mathcal L}_Y(D(\overline{a},0,0))\otimes_{{\mathcal O}_Y}{\mathcal O}_{W_S}\longrightarrow{\mathcal L}_Y(D(\overline{a},0,0))\otimes_{{\mathcal O}_Y}{\mathcal O}_{W}\longrightarrow0$$we see that it is enough to prove $$H^t(Y,{\mathcal L}_Y(D(\overline{a},0,0)-D_{M_0})\otimes_{{\mathcal O}_Y}{\mathcal O}_{W_S})=0$$for $t=1$. We do this for all $t>0$ by induction on $|S|$. The case $S=\emptyset$ was settled in Theorem \ref{wellvan}. If $S\ne\emptyset$ we pick an element $s\in S$, let $S'=S-\{s\}$ and consider the exact sequence$$0\longrightarrow{\mathcal L}_Y(D(\overline{a},0,0)-D_{M_0\cup\{s\}})\otimes_{{\mathcal O}_Y}{\mathcal O}_{W_{S'}}\longrightarrow{\mathcal L}_Y(D(\overline{a},0,0)-D_{M_0})\otimes_{{\mathcal O}_Y}{\mathcal O}_{W_{S'}}$$$$\longrightarrow{\mathcal L}_Y(D(\overline{a},0,0)-D_{M_0})\otimes_{{\mathcal O}_Y}{\mathcal O}_{W_{S}}\longrightarrow0.$$The induction hypothesis applies to the first two terms (note that also $M_0\cup\{s\}$ is stable !), hence also the third term has no higher cohomology.\hfill$\Box$\\

\section{The logarithmic de Rham complex on $Y$}
\label{secdera}

The results of this section are not needed in the remainder of this paper. The following lemma is an easy exercise (count the number of pairs $(ij)$ such that $a_{ij}$ for $(a_{ij})_{ij}\in {{{U}}}_{\tau}^{\sigma}$ is not yet forced to be zero by the requirements defining ${{{U}}}_{\tau}^{\sigma}$).

\begin{lem} Let $\tau=\{a_0,\ldots,a_r\}$ with $a_0<\ldots<a_r$, and let $\tau\subset\sigma$. Then $$|{{{U}}}_{\tau}^{\sigma}|=q^{\sum_{i=0}^ra_i-|\sigma^c\cap[0,a_i]|-i}.$$\hfill$\Box$\\
\end{lem}

For $\tau\subset\{1,\ldots,d\}$ we define $\overline{a}(\tau)=(\overline{a}_1(\tau),\ldots,\overline{a}_d(\tau))$ by\begin{gather}\overline{a}_i(\tau)=\left\{\begin{array}{l@{\quad:\quad}l}-1&i\in\tau\\0&i\in\{1,\ldots,d\}-\tau\end{array}\right.\label{atau}\end{gather}

\begin{lem}\label{dimdif} $$\dim_k\,H^0(Y,{\mathcal L}_Y(D(\overline{a}(\tau),0,0)))=q^{\sum_{i\in\tau}i}$$
\end{lem}

{\sc Proof:} First we claim$$\dim_k\,H^0(Y,{\mathcal L}_Y(D(\overline{a}(\tau),{{\bf{-1}}},0)))=\left\{\begin{array}{l@{\quad:\quad}l}0&\tau\ne\{1,\ldots,d\}\\\prod_{j=1}^d(q^j-1)&\tau=\{1,\ldots,d\}\end{array}\right.$$for ${{\bf{-1}}}=(-1,\ldots,-1)$. The argument for this is by outer induction on $d$ and inner induction on $\overline{a}_0(\tau)=-\sum_{j=1}^d\overline{a}_j(\tau)=|\tau|$. For $j\in \tau$ we have $$\dim_k\,H^0(Y,{\mathcal L}_Y(D(\overline{a}(\tau),{{\bf{-1}}},0)))=$$$$\dim_k\,H^0(Y,{\mathcal L}_Y(D(\overline{a}(\tau-\{j\}),{{\bf{-1}}},0)))+$$$$\sum_{j\in\sigma\in{\mathcal Y}}|{{{U}}}_{\sigma}-{{{U}}}_{\sigma}^{\{j\}}|\dim_k\,H^0(Y^{\sigma},{\mathcal L}_{Y^{\sigma}}(D(\overline{a}(\tau-\{j\})|_{\sigma},{{\bf{-1}}},0)))\dim_k\,H^0(Y^{\sigma^c},{\mathcal L}_{Y^{\sigma^c}}(D(\overline{a}(\tau)|_{\sigma^c},{{\bf{-1}}},0))).$$ Observing $\overline{a}(\tau-\{j\})=\overline{a}(\tau)^{[j]}$ and $\overline{a}(\tau-\{j\})|_{\sigma}=\overline{a}(\tau)|_{\sigma}^{[\iota_{\sigma}^{-1}(j)]}$ this follows from the proof of Theorem \ref{gallgvan}. By induction hypothesis this already settles the case $\tau\ne\{1,\ldots,d\}$. If $\tau=\{1,\ldots,d\}$ the easiest way is to choose $j=d$. In that case, by induction hypothesis, all summands on the right hand side except for $\sigma=\{d\}$ vanish, so we get $$\dim_k\,H^0(Y,{\mathcal L}_Y(D(\overline{a}(\tau),{{\bf{-1}}},0)))=$$$$|{{{U}}}_{\{d\}}-{{{U}}}_{\{d\}}^{\{d\}}|\dim_k\,H^0(Y^{\{0,\ldots,d-1\}},{\mathcal L}_{Y^{\{0,\ldots,d-1\}}}(D(\overline{a}(\tau)|_{\{0,\ldots,d-1\}},{{\bf{-1}}},0))).$$Since $|{{{U}}}_{\{d\}}-{{{U}}}_{\{d\}}^{\{d\}}|=q^d-1$ the induction hypothesis gives the claim. The claim established, the proof of the theorem itself is again by induction on $d$ and works precisely as the proof of Theorem \ref{wellvan}, with $S=\emptyset$ there. Namely, taking $W=V_{\{0\}}$, hence $T(W)=\{u.V_{\sigma};\,(\sigma,u)\in{\mathcal N}, 0\in\sigma\}$ so that $D_{T(W)}=D(0,{\bf 1},0)$, we saw that$$\dim_k\,H^0(Y,{\mathcal L}_Y(D(\overline{a}(\tau),0,0)))=$$$$\dim_k\,H^0(Y,{\mathcal L}_Y(D(\overline{a}(\tau),{{\bf{-1}}},0)))$$$$+\sum_{0\in\sigma}|\{u;\,(\sigma,u)\in{\mathcal N}\}|\dim_k\,H^0(Y^{\sigma},{\mathcal L}_{Y^{\sigma}}(D(\overline{a}(\tau)|_{\sigma},{{\bf{-1}}},0)))\dim_k\,H^0(Y^{\sigma^c},{\mathcal L}_{Y^{\sigma^c}}(D(\overline{a}(\tau)|_{\sigma^c},0,0))).$$All the terms are known by induction hypothesis, resp. the above claim, and we just have to sum up.\hfill$\Box$\\

Denote by $\Omega^{\bullet}_{Y}$ the de Rham complex on $Y$ with logarithmic poles along the normal crossings divisor $\sum_{V\in{\mathcal V}}V$ on $Y$. For $0\le r\le d$ we write $$z_r=\frac{\Xi_r}{\Xi_0},$$a rational function on $Y$. For $0\le s\le d$ denote by ${\mathcal P_s}$ the set of subsets of $\{1,\ldots,d\}$ consisting of $s$ elements. For a subset $\tau\subset\{1,\ldots,d\}$ let $$U(\tau)=\{(a_{ij})_{0\le i,j\le d}\in{{{U}}}(k)\quad|\quad a_{ij}=0\,\mbox{if}\, j\notin\{i\}\cup\tau \}.$$

\begin{satz}\label{logdif} (a) For each $0\le s\le d$ we have $H^t(Y,\Omega_Y^s)=0$ for all $t>0$.\\(b) The following set is a $k$-basis of $H^0(Y,\Omega_Y^s)$:
$$\{A.\bigwedge_{t\in\tau}\dlog(z_t)\quad|\quad \tau\in {\mathcal P_s}, A\in U(\tau)\}.$$In particular,$$\dim_k(H^0(Y,\Omega_Y^s))=\sum_{\tau\in{\mathcal P}_s}q^{\sum_{i\in\tau}i}.$$\end{satz}

{\sc Proof:} For elements $\tau^1=\{t^1_1<\ldots<t_s^1\}$ and $\tau^2=\{t^2_1<\ldots<t_s^2\}$ of ${\mathcal P_s}$ we write $\tau^1>\tau^2$ if and only if there is a $1\le b\le s$ such that $t^1_r=t^2_r$ for $1\le r\le b-1$ and $t^1_b>t^2_b$. In other words, we use the lexicographical ordering on ${\mathcal P_s}$. Let $$Y'=Y-\bigcup_{(\sigma,u)\in{\mathcal N}\atop u\ne1}u.V_{\sigma}.$$Writing $\Omega^{s}_{Y'}=\Omega^{s}_{Y}|_{Y'}$ we consider the following filtration of $\Omega^{s}_{Y'}$, indexed by ${\mathcal P_s}$:$$F^{\tau}\Omega^{s}_{Y'}=\bigoplus_{\tau'\in{\mathcal P_s}\atop \tau'\le\tau}\mathcal{O}_{Y'}\bigwedge_{t\in\tau'}dz_t.$$In particular $F^{\{d+1-s,\ldots,d\}}\Omega^{s}_{Y'}=\Omega^{s}_{Y'}$ if $s>0$. One checks that for open subsets $W_1$ and $W_2$ of $Y'$ and elements $g\in{{{U}}}(k)$ such that $gW_1=W_2$, the isomorphism $g:\Omega^{s}_{Y'}|_{W_1}\to\Omega^{s}_{Y'}|_{W_2}$ restricts to an isomorphism $F^{\tau}\Omega^{s}_{Y'}|_{W_1}\to F^{\tau}\Omega^{s}_{Y'}|_{W_2}$ for any $\tau$. Since ${{{U}}}(k)Y'=Y$ we may therefore extend this filtration from $Y'$ to all of $Y$, obtaining a ${{{U}}}(k)$-stable filtration$$(F^{\tau}\Omega^{s}_{Y})_{\tau\in{\mathcal P_s}}$$of $\Omega^{s}_{Y}$, with $F^{\{d+1-s,\ldots,d\}}\Omega^{s}_{Y}=\Omega^{s}_{Y}$ if $s>0$. We have\begin{align}Gr^{\tau}\Omega^{s}_{Y}&=\frac{F^{\tau}\Omega^{s}_{Y}}{\bigoplus_{\tau'<\tau}F^{\tau'}\Omega^{s}_{Y}}\cong{\mathcal L}_Y(D(\overline{a}(\tau),0,0))\notag \\{} &\quad\quad\quad f\bigwedge_{t\in\tau}dz_t\mapsto f.\notag\end{align}Hence we get $H^{t}(Y,Gr^{\tau}\Omega^{s}_{Y})=0$ for all $t>0$ from Corollory \ref{kritvan}. Since the $Gr^{\tau}\Omega^{s}_{Y}$ are isomorphic to the graded pieces of the filtration $(F^{\tau}\Omega^{s}_{Y})_{\tau\in{\mathcal P_s}}$ of $\Omega^{s}_{Y}$ this implies statement (a). Moreover, the classes of the elements $A.\bigwedge_{t\in\tau}\dlog(z_t)$ for $A\in U(\tau)$ form a $k$-basis of $H^0(Y,Gr^{\tau}\Omega^{s}_{Y})$. This can be proven by tracing back the proof of Lemma \ref{dimdif}; or by showing the linear independence of the $A.\bigwedge_{t\in\tau}\dlog(z_t)$ and {\it applying} Lemma \ref{dimdif}. We get statement (b).\hfill$\Box$\\

\begin{kor}\label{hodgedeg} The Hodge spectral sequence $E_1^{st}=H^t(Y,\Omega_Y^s)\Longrightarrow H^{s+t}(Y,\Omega_Y^{\bullet})$ degenerates in $E_1$. We have $H^s(Y,\Omega^{\bullet}_Y)=H^0(Y,\Omega^{s}_Y)$ for all $s$.
\end{kor}

We wish to determine $H^0(Y,\Omega_Y^s)$ as a ${\rm GL}_{d+1}(k)$-representation. For this we recall the classification of irreducibel representations of ${\rm GL}_{d+1}(k)$ on ${k}$-vector spaces according to Carter and Lusztig. For $1\le r\le d$ let $t_r\in{\rm GL}_{d+1}(k)$ denote the permutation matrix obtained by interchanging the $(r-1)$-st and the $r$-th row (or equivalently: column) of the identity matrix (recall that we start counting with $0$). Then $S=\{t_1,\ldots,t_d\}$ is a set of Coxeter generators for the Weyl group of ${\rm GL}_{d+1}(k)$. 

\begin{satz} \cite{calu} (i) For an irreducible representation $\rho$ of ${\rm GL}_{d+1}(k)$ on a ${k}$-vector space, the subspace $\rho^{U(k)}$ of $U(k)$-invariants is one dimensional. If the action of the group $B(k)$ of upper triangular matrices on $\rho^{U(k)}$ is given by the character $\chi:B(k)/U(k)\to{k}^{\times}$ and if $J=\{t\in S;\,t.\rho^{U(k)}=\rho^{U(k)}\}$, then the pair $(\chi,J)$ determines $\rho$ up to isomorphism.\\(ii) Conversely, given a character $\chi:B(k)/U(k)\to{k}^{\times}$ and a subset $J$ of $\{t\in S;\,\chi^t=\chi\}$, there exists an irreducible representation $\Theta(\chi,J)$ of ${\rm GL}_{d+1}(k)$ on a ${k}$-vector space whose associated pair (as above) is $(\chi,J)$. 
\end{satz}

For $1\le j\le d$ we need the rational function \begin{gather}\gamma_j=\prod_{(a_0,\ldots,a_{j-1})\in k^j}(z_j+a_{j-1}z_{j-1}+\ldots+a_1z_1+a_0)\label{gammaj}\end{gather} on $Y$, and if in addition $0\le s\le d$ we define the integer$$m_j^s=\max\{0,s-j\}q -\max\{0,s-j+1\}.$$

\begin{satz}\label{unimod} For $0\le s\le d$, the ${\rm GL}_{d+1}(k)$-representation on $H^0(Y,\Omega_Y^s)$ is equivalent to $\Theta(1,\{t_{s+1},\ldots,t_d\})$. The subspace of $U(k)$-invariants of $H^0(Y,\Omega_Y^s)$ is generated by $$\omega_s=(\prod_{j=1}^d\gamma_j^{m_j^s})dz_1\wedge\ldots\wedge dz_s.$$
\end{satz}

{\sc Proof:} (i) We check that $\omega_s$ is indeed an element of $H^0(Y,\Omega_Y^s)$. In notations of section \ref{cohmodp} let us abbreviate$$W_{\sigma}=\sum_{u\in U_{\sigma}}u.V_{\sigma}$$for $\sigma\in{\mathcal{Y}}$, a reduced and $U(k)$-stable divisor. For $1\le j\le d$ one finds that the zero-pole divisor of $\gamma_j^{-1}$ is\begin{gather}\sum_{\sigma\atop j\in\sigma, 0\notin\sigma}(q^{|\sigma\cap[0,j-1]|})W_{\sigma}+\sum_{\sigma\atop j\in\sigma,0\in\sigma}(q^{|\sigma\cap[0,j-1]|}-q^j)W_{\sigma}+\sum_{\sigma\atop j\notin\sigma,0\in\sigma}-q^jW_{\sigma}.\label{gadiv}\end{gather}Hence the zero pole divisor of $\prod_{j=1}^d\gamma_j^{m_j^s}$ is$$\sum_{\sigma\atop0\notin\sigma}\sum_{j\in\sigma}-m_j^sq^{|\sigma\cap[0,j-1]|}W_{\sigma}+\sum_{\sigma\atop0\in\sigma}(\sum_{j\ge1,j\in\sigma}-m_j^sq^{|\sigma\cap[0,j-1]|}+\sum_{j=1}^dm_j^sq^j)W_{\sigma}.$$This divisor is smaller (for the usual order on the set of Cartier divisors) than$$D_s=\sum_{\sigma\atop 0\notin\sigma}|\sigma\cap[1,s]|W_{\sigma}+\sum_{\sigma\atop 0\in\sigma}(|\sigma\cap[1,s]-s)W_{\sigma}.$$Now we have the ($U(k)$-equivariant) embedding$${\mathcal L}_Y(D_s)\longrightarrow\Omega_Y^s,\quad f\mapsto fdz_1\wedge\ldots\wedge dz_s.$$By the above, $\prod_{j=1}^d\gamma_j^{m_j^s}$ is a global section of ${\mathcal L}_Y(D_s)$; its image in $H^0(Y,\Omega_Y^s)$ is $\omega_s$.\\(ii) We check that $\omega_s$ is fixed by $B(k)$. The rational function $\gamma_j$ is $U(k)$-stable, and similarly the $s$-form $dz_1\wedge\ldots\wedge dz_s$ is $U(k)$-stable. On the other hand, the torus $T(k)$ of diagonal matrices in ${\rm GL}_{d+1}(k)$ acts as follows: if $t\in T(k)$ has diagonal entries $a_{00},\ldots,a_{dd}$, then $t$ acts on $\gamma_j$ by permuting its factors different from $z_j$, and on the factor $z_j$ it acts by multiplication with $a_{jj}/a_{00}$. Hence $t.\gamma_j^{m_j^s}=(a_{00}/a_{jj})\gamma_j^{m_j^s}$ --- because $m_j^s\equiv -1$ (modulo $q$) --- and $t.dz_j=(a_{jj}/a_{00})dz_j$. Multiplying together we get our claim.\\(iii) It is clear that $\{t_{s+1},\ldots,t_d\}=\{t\in S;\,t.\omega_s=\omega_s\}$.\\(iv) From (i), (ii), (iii) it follows that $H^0(Y,\Omega_Y^s)$ contains $\Theta(1,\{t_{s+1},\ldots,t_d\})$. We have$$\dim_k(\Theta(1,\emptyset))=q^{1+2+\ldots+d}=\dim_k(H^0(Y,\Omega_Y^d))$$where the first equality is \cite{caen} Theorem 6.12(ii) and the second one follows from Theorem \ref{logdif}; hence our statement in the case $s=d$. To go further we point out the following trivial consequence: the $k[{\rm GL}_{d+1}(k)]$-module generated by $\omega_d$ contains a non-zero logarithmic differential $d$-form (because $H^0(Y,\Omega_Y^d)$ contains such a form, by Theorem \ref{logdif}). For $s<d$ we know from Theorem \ref{logdif} that $H^0(Y,\Omega_Y^s)$ is generated as a $k$-vector space by logarithmic differential $s$-forms. Since ${\rm GL}_{d+1}(k)$ acts transitively on the set of $k^{\times}$-homothety classes of non-zero logarithmic differential $s$-forms, each non-zero logarithmic differential $s$-form generates $H^0(Y,\Omega_Y^s)$ as a $k[{\rm GL}_{d+1}(k)]$-module. Hence we need to show that the sub-$k[{\rm GL}_{d+1}(k)]$-module generated by $\omega_s$ contains a non-zero logarithmic differential $s$-form. As a formal expression, $\omega_s$ is independent of $d$, as long as $d\ge s$; this follows immediately from the definition of the numbers $m_j^s$. Moreover, the action of the subgroup ${\rm GL}_{s+1}(k)$ on $\omega_s$ is independent of $d$ if we use the obvious embedding ${\rm GL}_{s+1}(k)\to {\rm GL}_{d+1}(k)$ (into the upper left square, followed by entries $=1$ on the rest of the diagonal). Therefore the consequence pointed out above (when $d$ there takes the value of our present $s$) says that the $k[{\rm GL}_{s+1}(k)]$-module generated by $\omega_s$ contains a non-zero logarithmic differential $s$-form; in particular, the $k[{\rm GL}_{d+1}(k)]$-module generated by $\omega_s$ contains a non-zero logarithmic differential $s$-form. We are done.\hfill$\Box$\\

We draw representation theoretic consequences. Fix $0\le s\le d$, let $n_1=\ldots=n_s=1$ and $n_{s+1}=d+1-s$ and define the parabolic subgroup$$P_{s}=\{(a_{ij})_{0\le i,j\le d}\in{\rm GL}_{d+1};\quad a_{ij}=0\mbox{ if }j+1\le n_1+\ldots+n_{\ell}\mbox{ and }i+1>n_1+\ldots+n_{\ell}\mbox{ for some }\ell\}$$of ${\rm GL}_{d+1}$. For any group $H$ let ${\bf 1}$ denote the trivial $k[H]$-module.

\begin{kor} The following generalized Steinberg representation is irreducible:$${\rm ind}_{P_{s}(k)}^{{\rm GL}_{d+1}(k)}{\bf 1}/\sum_{P\supsetneq P_{s}}{\rm ind}_{P(k)}^{{\rm GL}_{d+1}(k)}{\bf 1}$$where the sum runs over all parabolic subgroups $P$ of ${\rm GL}_{d+1}$ strictly containing $P_{s}$. 
\end{kor}

{\sc Proof:} This representation clearly contains $\Theta(1,\{t_{s+1},\ldots,t_d\})$ and its $k$-dimension is $\sum_{\tau\in{\mathcal P}_s}q^{\sum_{i\in\tau}i}$ (to see this use for example the formula \cite{caen} exc. 6.3). But this is also the $k$-dimension of $\Theta(1,\{t_{s+1},\ldots,t_d\})$ as follows from Theorem \ref{unimod}.\hfill$\Box$\\

Let $W$ denote the ring of Witt vectors with coefficients in $k$. Endow $\spec(k)$ and $\spf(W)$ with the trivial log structure and endow $Y$ with the log structure corresponding to the normal crossings divisor $\cup_{V\in{\mathcal V}}V$ on $Y$. We ask for the logarithmic crystalline cohomology $H_{crys}^*(Y/W)$ of $Y/k$ relative to the divided power thickening $\spf(W)$ of $\spec(k)$. (Note that $H_{crys}^*(Y/W)\otimes{\mathbb Q}$ is the rigid cohomology of the open subscheme $Y-\cup_{V\in{\mathcal V}}V$ of $Y$.)

\begin{satz}\label{bacauni} $H_{crys}^s(Y/W)$ is torsion free for any $s$, and$$H_{crys}^s(Y/W)\otimes_{W}k=H^s(Y,\Omega_Y^{\bullet})=H^0(Y,\Omega_Y^{s}).$$
\end{satz}

{\sc Proof:} We already remarked $H^s(Y,\Omega_Y^{\bullet})=H^0(Y,\Omega_Y^{s})$. The base change property $H_{crys}^s(Y/W)\otimes_{W}k=H^s(Y,\Omega_Y^{\bullet})$ is a consequence of the torsion freeness of $H_{crys}^s(Y/W)$ (for all $s$). To prove it it suffices by \cite{blka} (7.3) to prove that the log scheme $Y$ is ordinary, i.e. that the Newton polygon of $H_{crys}^s(Y/W)\otimes{\mathbb Q}$ coincides with the Hodge polygon of $H^s(Y,\Omega_Y^{\bullet})$. But the endpoints of these two polygons are the same, and both have only a single slope: for the Hodge polygon this is our equality $H^s(Y,\Omega_Y^{\bullet})=H^0(Y,\Omega_Y^{s})$, for the Newton polygon this was verified in \cite{phien} Theorem 6.3.\hfill$\Box$\\ 

{\it Remarks:} (a) See sections \ref{linsec} and \ref{hosec} for the notations in this remark which have not yet been defined. Let $\Omega_{\mathfrak X}^s$ be the degree $s$ term of the relative logarithmic de Rham complex of ${\mathfrak X}$ over $\spf({\mathcal O}_K)$ (with respect to the log structures defined by the special fibres). Then $\Omega_{\mathfrak X}^s\otimes_{{\mathcal O}_K}{\mathcal O}_{\widehat{K}}\cong {\mathcal V}_{{\mathcal O}_{\widehat{K}}}$ for the irreducible rational representation $V$ of $L_1$ with heighest weight$$\sum_{i=1}^s(s+1)\epsilon_i+\sum_{i=s+1}^ds\epsilon_i\in X^*(T_1).$$(In fact the remark following the proof of Proposition \ref{intgstab} applies and we have $\Omega_{\mathfrak X}^s\cong {\mathcal V}_{{\mathcal O}_{{K}}}$). On the other hand we have $\Omega_{\mathfrak X}^s\otimes_{{\mathcal O}_K}{\mathcal O}_Y=\Omega_{Y}^s$. Via these isomorphisms, the filtration appearing in the proof of Theorem \ref{logdif} is just the one appearing in the proof of \ref{globnull}, and statement (a) in Theorem \ref{logdif} is statement (\ref{grsimpvan}) in the proof of \ref{globnull} (for this $V$).

(b) Theorem \ref{logdif} is used in \cite{latt} to prove some new cases of Schneider's conjecture (see \cite{schn}) concerning $p$-adic analytic splittings of Hodge filtrations of the de Rham cohomology of certain local systems on projective varieties uniformized by $X$. 

\section{Lattices in line bundles on the symmetric space}

\label{linsec}

Let $T$ be the torus of diagonal matrices in $G$ and let $X_*(T)$, resp. $X^*(T)$, denote the group of algebraic cocharacters, resp. characters, of $T$. Put $A=X_*(T)\otimes{\mathbb{R}}.$ For $0\le i\le d$ define the cocharacters$$e_i:{\mathbb{G}}_m\to{\rm GL}\sb {d+1},\quad t\mapsto\diag(1,\ldots,1,t,1,\ldots,1)$$with $t$ as the $i$-th diagonal entry, i.e. $e_i(t)_{ii}=t$, $e_i(t)_{jj}=1$ for $i\ne j$ and $e_i(t)_{j_1j_2}=0$ for $j_1\ne j_2$. The $e_i$ form a basis of the $\mathbb{R}$-vector space $A$, hence an identification $A=\mathbb{R}^{d+1}$. The pairing $X_*(T)\times X^*(T)\to\mathbb{Z}$ which sends $(x,\mu)$ to the integer $\mu(x)$ such that $\mu(x(y))=y^{\mu(x)}$ for any $y\in\mathbb{G}_m$ extends to a duality between $\mathbb{R}$-vector spaces$$A\times(X^*(T)\otimes\mathbb{R})\longrightarrow\mathbb{R}$$$$(x,\mu)\mapsto\mu(x).$$Let $\epsilon_0,\ldots,\epsilon_{d}\in X^*(T)$ denote the basis dual to $e_0,\ldots,e_{d}$.

Let ${\mathbb A}((K^{d+1})^*)$, resp. ${\mathbb P}((K^{d+1})^*)$, denote the affine, resp. the projective space spanned by $(K^{d+1})^*=\Hom_K(K^{d+1},K)$. The action of $G={\rm GL}_{d+1}(K)={\rm GL}(K^{d+1})$ on $(K^{d+1})^*$ defines an action of $G$ on the Drinfel'd symmetric spaces$$X^{cone}={\mathbb A}((K^{d+1})^*)-(\mbox{the union of all $K$-rational hyperplanes through}\,0),$$$$X={\mathbb P}((K^{d+1})^*)-(\mbox{the union of all $K$-rational hyperplanes}).$$There is a natural $G$-equivariant projection of $K$-rigid spaces $X^{cone}\to X$. Let $\Xi_0,\ldots,\Xi_{d}$ be the standard coordinate functions on $X^{cone}$ corresponding to the canonical basis of $(K^{d+1})^*$; they induce a set of projective coordinate functions on $X$.

Let ${\mathfrak{X}}$ be the strictly semistable formal ${\mathcal O}_K$-scheme with generic fibre $X$ introduced in \cite{mus}. The action of $G$ on $X$ extends naturally to ${\mathfrak{X}}$. For open formal subschemes ${\mathfrak{U}}$ of ${\mathfrak{X}}$ we put$${\mathcal O}_{\widehat{{\mathfrak{U}}}}={\mathcal O}_{\mathfrak{U}}\otimes_{{\mathcal O}_K}{\mathcal O}_{\widehat{K}}.$$For $j\ge 0$ let $F^j$ be the set of non-empty intersections of $(j+1)$-many pairwise distinct irreducible components of ${\mathfrak{X}}\otimes_{{\mathcal O}_K}k$. (Thus $F^j$ is in natural bijection with the set of $j$-simplices of the Bruhat Tits building of ${\rm PGL}\sb {d+1}/K$.) For $Z\in F^0$ let ${\mathfrak{U}}_Z$ be the maximal open formal subscheme of ${\mathfrak X}$ such that ${\mathfrak{U}}_Z\otimes_{{\mathcal O}_K}k$ is contained in $Z$. Let $Y\in F^0$ be the central irreducible component of ${\mathfrak X}\otimes_{{\mathcal O}_K}k$ with respect to $\Xi_0,\ldots,\Xi_d$, characterized by the following condition. By construction we may view $\Xi_i\Xi_0^{-1}$ for $0\le i\le d$ as a global section of ${\mathcal O}_{\mathfrak X}\otimes_{{\mathcal O}_K}K$ (i.e. as a rigid analytic function on $X$). Now $Y\in F^0$ is the unique irreducible component such that for all unimodular $d+1$-tupel $(a_0,\ldots,a_d)\in{\mathcal O}_K^{d+1}$ (at least one $a_i$ is a unit in ${\mathcal O}_K$) the linear combination $\sum_{i=0}^da_i\Xi_i\Xi_0^{-1}$, when restricted to ${\mathfrak{U}}_Y$, in fact lies in the subalgebra ${\mathcal O}_{\mathfrak X}({\mathfrak{U}}_Y)$ of ${\mathcal O}_{\mathfrak X}\otimes_{{\mathcal O}_K}K({\mathfrak{U}}_Y)$ and is even invertible in ${\mathcal O}_{\mathfrak X}({\mathfrak{U}}_Y)$. By \cite{mus} we may identify this $k$-scheme $Y$ with the one from section \ref{cohmodp}. The subgroup $K^{\times}.{\rm GL}\sb {d+1}({{\mathcal O}_K})$ of $G$ is the stabilizer of $Y$. We define the subset $$F^0_{A}=T.Y,$$the orbit of $Y$ for the action of $T$ on the set $F^0$. (This set corresponds to the set of vertices in an apartment of the Bruhat Tits building of ${\rm PGL}\sb {d+1}/K$.) We denote by ${\mathfrak Y}$ the maximal open formal subscheme of ${\mathfrak X}$ such that ${\mathfrak Y}\otimes k$ is contained in the closed subscheme $\cup_{Z\in F^0_A}Z$ of ${\mathfrak X}\otimes k$. Note that the open subscheme$$Y'={\mathfrak Y}\cap Y$$of $Y$ is the complement of all divisors $u.V_{\sigma}$ for $(\tau,u)\in{\mathcal{N}}$ with $u\ne1$. Also observe that $U(K).{\mathfrak Y}={\mathfrak X}$ and $U(k).Y'=Y$. Let $\overline{T}=T/K^{\times}$. For $\mu=\sum_{j=0}^da_j\epsilon_j\in X^*(T)$ let \begin{gather}\overline{\mu}=(\frac{1}{d+1}\sum_{j=0}^da_j)(\sum_{j=0}^d\epsilon_j)-\mu,\label{mubar}\end{gather}an element of the subspace $X^*(\overline{T})\otimes\mathbb{R}$ of $X^*({T})\otimes\mathbb{R}$. Letting \begin{gather}\overline{a}_{j}(\mu)=\frac{(\sum_{i\ne j}a_i)-da_{j}}{d+1}\label{aquer}\end{gather}we have$$\overline{\mu}=\sum_{j=0}^d\overline{a}_{j}(\mu)\epsilon_j.$$
There is a $\delta(\mu)\in\frac{1}{d+1}.\mathbb{Z}\cap[0,1[$ such that $\overline{a}_{j}(\mu)+\delta(\mu)\in\mathbb{Z}$ for all $0\le j\le d$. If $\delta(\mu)=0$ --- this is equivalent with $\overline{\mu}\in X^*(\overline{T})$ --- we define $\overline{a}(\mu)=(\overline{a}(\mu)_i)_{1\le i\le d}\in{\mathbb{Z}}^d$ by setting $$\overline{a}(\mu)_i=\overline{a}_i(\mu)$$and then let, as in section \ref{cohmodp}, $\overline{a}(\mu)_0=-\sum_{i=1}^d\overline{a}(\mu)_i=\overline{a}_0(\mu)$. If $\delta(\mu)\ne0$ --- this is equivalent with $\overline{\mu}\notin X^*(\overline{T})$ ---  we define $\lceil\overline{a}(\mu)\rceil=(\lceil\overline{a}(\mu)\rceil_i)_{1\le i\le d}\in{\mathbb{Z}}^d$ by setting $$\lceil\overline{a}(\mu)\rceil_i=\lceil\overline{a}_i(\mu)\rceil$$and then let, as in section \ref{cohmodp}, $\lceil\overline{a}(\mu)\rceil_0=-\sum_{i=1}^d\lceil\overline{a}(\mu)\rceil_i=-d\delta(\mu)-\sum_{i=1}^d\overline{a}_i(\mu)$.

We define ${m}(\mu)=(m(\mu)_i)_{1\le i\le d}\in {\mathbb{Z}}^d$ and ${n}(\mu)=(n(\mu)_i)_{1\le i\le d}\in {\mathbb{Z}}^d$ by setting $$m(\mu)_i=\lfloor i\delta(\mu)\rfloor\quad\quad\mbox{ and }\quad\quad n(\mu)_i=\lfloor(i-1-d)\delta(\mu)\rfloor.$$For $Z\in F^0_{A}$ and $\gamma\in X^*(\overline{T})\otimes{\mathbb{R}}$ we define $\gamma(Z)\in\mathbb{R}$ as$$\gamma(Z)=\gamma(t)\quad\mbox{with}\,\,t\in T\,\,\mbox{such that}\,\,t.Y=Z.$$For $Z\in F^0_{A}$ let ${\mathcal J}_{Z}\subset {\mathcal O}_{\mathfrak Y}$ be the ideal defining $Z\cap{\mathfrak Y}$ inside ${\mathfrak Y}$. Note that ${\mathcal J}_{Z}$ is invertible inside ${\mathcal O}_{\mathfrak Y}\otimes_{{\mathcal O}_K}K$: indeed, small open formal subschemes of ${\mathfrak Y}$ admit open embeddings into the $\pi$-adic completion of $\spec({\mathcal O}_K[X_0,\ldots,X_d]/(X_0\ldots X_d-\pi))$, and there a typical generator of ${\mathcal J}_{Z}$ is of the form $X_0$ (for an appropriate numbering of $X_0,\ldots,X_d$); in $K[X_0,\ldots,X_d]/(X_0\ldots X_d-\pi)$ its inverse is $\pi^{-1}X_1\ldots X_d$. Thus we may speak of negative integral powers of ${\mathcal J}_{Z}$ as ${\mathcal O}_{\mathfrak Y}$-submodules of ${\mathcal O}_{\mathfrak Y}\otimes_{{\mathcal O}_K}K$. Also note that on small open formal subschemes of ${\mathfrak Y}$ we have ${\mathcal J}_{Z}={\mathcal O}_{\mathfrak Y}$ for almost all $Z$, therefore the infinite products of ${\mathcal O}_{\widehat{\mathfrak Y}}$-submodules inside ${\mathcal O}_{\mathfrak Y}\otimes_{{\mathcal O}_K}{\widehat{K}}$ below make sense. On ${\mathfrak Y}$ we define the subsheaf\begin{gather}({\mathcal O}_{\widehat{\mathfrak Y}})^{\overline{\mu}}=\sum_{s=0}^d\widehat{\pi}^s\prod_{Z\in F^0_A}{\mathcal J}_{Z}^{\lceil\overline{\mu}(Z)-\frac{s}{d+1}\rceil},\label{iii}\end{gather}of ${\mathcal O}_{{\mathfrak{Y}}}\otimes_{{\mathcal O}_K}{\widehat{K}}$, i.e. $({\mathcal O}_{\widehat{\mathfrak Y}})^{\overline{\mu}}$ is the ${\mathcal O}_{\widehat{\mathfrak Y}}$-submodule of ${\mathcal O}_{{\mathfrak{Y}}}\otimes_{{\mathcal O}_K}{\widehat{K}}$ generated by the ${\mathcal O}_{\widehat{\mathfrak Y}}$-submodules $$\widehat{\pi}^s\prod_{Z\in F^0_{A}}{\mathcal J}_{Z}^{\lceil\overline{\mu}(Z)-\frac{s}{d+1}\rceil}\quad\quad(s=0,\ldots,d).$$Let $({\mathcal O}_{\widehat{\mathfrak X}})^{\overline{\mu}}$ be the unique ${U}(K)$-equivariant subsheaf of ${\mathcal O}_{{\mathfrak{X}}}\otimes_{{\mathcal O}_K}{\widehat{K}}$ (with its ${U}(K)$-action induced by that of ${U}(K)\subset G$ on ${\mathfrak{X}}$) whose restriction to ${\mathfrak Y}$ is $({\mathcal O}_{\widehat{\mathfrak Y}})^{\overline{\mu}}$.

\begin{lem}\label{lbdlext} (a) If $\overline{\mu}\in X^*(\overline{T})$ then $({\mathcal O}_{\widehat{\mathfrak X}})^{\overline{\mu}}$ is a line bundle on ${\mathfrak X}$. We have an isomorphism\begin{gather}{\mathcal L}_{Y}(D(\overline{a}(\mu),0,0))\cong ({\mathcal O}_{\widehat{\mathfrak X}})^{\overline{\mu}}\otimes_{{\mathcal O}_{\widehat{\mathfrak{X}}}}{\mathcal O}_{Y}.\label{lbdldenu}\end{gather}
(b) If $\overline{\mu}\notin X^*(\overline{T})$ we have isomorphisms\begin{gather}({\mathcal O}_{\widehat{\mathfrak X}})^{\overline{\mu}}\otimes_{{\mathcal O}_{\widehat{K}}}k\cong\prod_{x\in F^0}\frac{({\mathcal O}_{\widehat{\mathfrak X}})^{\overline{\mu}}\otimes_{{\mathcal O}_{\widehat{\mathfrak{X}}}}{\mathcal O}_{Z}}{{\mathcal O}_{Z}\mbox{-torsion}},\label{lbdlzerldeninu}\end{gather}\begin{gather}{\mathcal L}_{Y}(D(\lceil\overline{a}(\mu)\rceil,n(\mu),m(\mu)))\cong \frac{({\mathcal O}_{\widehat{\mathfrak X}})^{\overline{\mu}}\otimes_{{\mathcal O}_{\widehat{\mathfrak{X}}}}{\mathcal O}_{Y}}{{\mathcal O}_{Y}\mbox{-torsion}}.\label{lbdldeninu}\end{gather}
\end{lem}

{\sc Proof:} (a) Here $\overline{\mu}(Z)\in\mathbb{Z}$ for all $Z\in F^0_A$ and it follows from formula (\ref{iii}) that\begin{gather}({\mathcal O}_{\widehat{\mathfrak Y}})^{\overline{\mu}}=\prod_{Z\in F^0_A}{\mathcal J}_{Z}^{\overline{\mu}(Z)}.\label{deltanu}\end{gather}In particular, $({\mathcal O}_{\widehat{\mathfrak Y}})^{\overline{\mu}}$ is a line bundle on ${\mathfrak Y}$ in that case. Now $({\mathcal O}_{\widehat{\mathfrak X}})^{\overline{\mu}}\otimes_{{\mathcal O}_{\widehat{\mathfrak{X}}}}{\mathcal O}_{Y}$ with its ${{{U}}}(k)$-action is the unique ${{{U}}}(k)$-equivariant subsheaf of the constant sheaf with value the function field $k(Y)$ on $Y$ whose restriction to $Y'$ is $$({\mathcal O}_{\widehat{\mathfrak Y}})^{\overline{\mu}}\otimes_{{\mathcal O}_{\widehat{\mathfrak{Y}}}}{\mathcal O}_{Y'}$$(for the uniqueness note that ${{{U}}}(k).Y'=Y$). Write $({\mathcal O}_{\widehat{\mathfrak Y}})^{\overline{\mu}}\otimes_{{\mathcal O}_{\widehat{\mathfrak{Y}}}}{\mathcal O}_{Y'}={\mathcal L}_{Y'}(D)$ (as subsheaves of the constant sheaf $k(Y)$ on $Y'$) with a divisor $D$ on $Y'$. By ${{{U}}}(k)$-equivariance of its both sides and ${{{U}}}(k).Y'=Y$, to prove formula (\ref{lbdldenu}) we only need to prove $D=D(\overline{a}(\mu),0,0)|_{Y'}$. This holds because for $\emptyset\ne \sigma\subsetneq\Upsilon=\{0,\ldots,d\}$ the divisor $V_{\sigma}$ on $Y$ as defined in section \ref{cohmodp} is the divisor $Z_{{\sigma}}\cap Y$ on $Y$. Here we write $Z_{{\sigma}}= t_{\sigma}Y\in F^0_A$ with $t_{\sigma}=\diag(t_{\sigma,0},\ldots,t_{\sigma,d})\in T\subset G$ defined as $t_{\sigma,j}=1$ if $j\notin\sigma$ and $t_{\sigma,j}=\pi$ if $j\in\sigma$. We then apply equation (\ref{deltanu}) which tells us that the prime divisor $V_{\sigma}\cap Y'$ occurs with multiplicity $-\overline{\mu}(Z_{\sigma})=-\sum_{j\in\sigma}\overline{a}(\mu)_j$ in $D$.

(b) We proceed similarly as in case (a). For $Z\in F^0_{A}$ define the number $s(Z)\in[0,d]\cap\mathbb{Z}$ by requiring $\overline{\mu}(Z)-\frac{s(Z)}{d+1}\in\mathbb{Z}$. The definition of $\overline{\mu}$ together with our assumption $\overline{\mu}\notin X^*(\overline{T})$ implies $s(Z)\ne s(Z')$ for any two distinct but neighbouring $Z, Z'$ in $F^0_A$. Since restricted to ${\mathfrak U}_Z$ the relevant summand in formula (\ref{iii}) is the one for $s=s(Z)$ we obtain: the reduction $({\mathcal O}_{\widehat{\mathfrak Y}})^{\overline{\mu}}\otimes_{{\mathcal O}_{\widehat{K}}}k$ of $({\mathcal O}_{\widehat{\mathfrak Y}})^{\overline{\mu}}$ decomposes into a product, indexed by the set $F^0_{A}$ whose factor for $Z\in F^0_{A}$ is the image of the map$$\widehat{\pi}^{s(Z)}\prod_{Z'\in F^0_{A}}{\mathcal J}_{Z'}^{\lceil\overline{\mu}(Z')-\frac{s(Z)}{d+1}\rceil}\longrightarrow({\mathcal O}_{\widehat{\mathfrak Y}})^{\overline{\mu}}\longrightarrow{({\mathcal O}_{\widehat{\mathfrak Y}})^{\overline{\mu}}}\otimes_{{\mathcal O}_{\widehat{\mathfrak{Y}}}}{\mathcal{O}}_{Z\cap {\mathfrak Y}}.$$This is a line bundle on $Z\cap {\mathfrak Y}$ and maps isomorphically to the quotient of ${({\mathcal O}_{\widehat{\mathfrak Y}})^{\overline{\mu}}}\otimes_{{\mathcal O}_{\widehat{\mathfrak{Y}}}}{\mathcal{O}}_{Z\cap {\mathfrak Y}}$ divided by its ${\mathcal O}_{Z\cap {\mathfrak Y}}$-torsion. Thus\begin{gather}({\mathcal O}_{\widehat{\mathfrak Y}})^{\overline{\mu}}\otimes_{{\mathcal O}_{\widehat{K}}}k\cong\prod_{Z\in F^0_A}\frac{({\mathcal O}_{\widehat{\mathfrak Y}})^{\overline{\mu}}\otimes_{{\mathcal O}_{\widehat{\mathfrak{Y}}}}{\mathcal{O}}_{Z\cap {\mathfrak Y}}}{{\mathcal{O}}_{Z\cap {\mathfrak Y}}\mbox{-torsion}}\label{deltaninu}\end{gather}and each $$\frac{({\mathcal O}_{\widehat{\mathfrak Y}})^{\overline{\mu}}\otimes_{{\mathcal O}_{\mathfrak{Y}}}{\mathcal{O}}_{Z\cap {\mathfrak Y}}}{{\mathcal{O}}_{Z\cap {\mathfrak Y}}\mbox{-torsion}}$$is an invertible ${\mathcal{O}}_{Z\cap {\mathfrak Y}}$-module. Hence formula (\ref{lbdlzerldeninu}). If we define the divisor $D$ on $Y'$ by requiring $$\frac{({\mathcal O}_{\widehat{\mathfrak Y}})^{\overline{\mu}}\otimes_{{\mathcal O}_{\widehat{\mathfrak{Y}}}}{\mathcal{O}}_{Y'}}{{\mathcal{O}}_{Y'}\mbox{-torsion}}={\mathcal L}_{Y'}(D)$$then to prove formula (\ref{lbdldeninu}) we need to prove $D=D(\lceil\overline{a}(\mu)\rceil,n(\mu),m(\mu))|_{Y'}$. This holds because for $\emptyset\ne\sigma\subset\Upsilon-\{0\}$ the prime divisor $V_{\sigma}$ occurs in $D$ with multiplicity $$-\lceil\overline{\mu}(Z_{\sigma})\rceil=-\lceil\sum_{j\in\sigma}\overline{a}_j(\mu)\rceil=m(\mu)_{|\sigma|}-\sum_{j\in\sigma}\lceil\overline{a}(\mu)\rceil_j$$and for $\sigma\subsetneq\Upsilon$ with $0\in\sigma$ the prime divisor $V_{\sigma}$ occurs in $D$ with multiplicity $$-\lceil\overline{\mu}(Z_{\sigma})\rceil=-\lceil\sum_{j\in\sigma}\overline{a}_j(\mu)\rceil=n(\mu)_{|\sigma|}-\sum_{j\in\sigma}\lceil\overline{a}(\mu)\rceil_j$$with $Z_{\sigma}$ as defined before.\hfill$\Box$\\

\section{The holomorphic discrete series}
\label{hosec}

Let $$\Phi=\{\epsilon_i-\epsilon_j;\,\,0\le i,j\le d\,\,\mbox{and}\,i\ne j\}\subset X^*(\overline{T}).$$For $0\le i,j\le d$ and $i\ne j$ define the morphism of algebraic groups over $\mathbb{Z}$\begin{gather}\widetilde{\alpha}_{ij}:\mathbb{G}_a\longrightarrow {\rm GL}\sb {d+1},\quad u\mapsto I_{d+1}+u.e_{ij}\label{alphaij}\end{gather}where $I_{d+1}+u.e_{ij}$ is the matrix $(u_{rs})$ with $u_{rr}=1$ (all $r$), with $u_{ij}=u$ and with $u_{rs}=0$ for all other pairs $(r,s)$. For the root $\alpha=\epsilon_i-\epsilon_j\in \Phi$ and $r\in\mathbb{R}$ let $$U_{\alpha,r}=\widetilde{\alpha}_{ij}(\{u\in K;\,\omega(u)\ge r\})\subset G.$$ For $x\in A$ let$$U_x=\mbox{the subgroup of}\,G\,\mbox{generated by all}\,U_{\alpha,-\alpha(x)}\,\mbox{for}\,\alpha\in\Phi.$$We recall definitions from \cite{schtei}. Note that many conventions are opposite to those in \cite{schn}. Define the ${\rm GL}\sb {d+1}$-subgroups$$P_1=\mbox{matrices of the form}\quad\left(\begin{array}{cccc}1&0&\cdots&0\\{*}&*&\cdots&*\\\vdots&\vdots&{}&\vdots\\{*}&*&\cdots&*\end{array}\right),$$$$L_1=\mbox{matrices of the form}\quad\left(\begin{array}{cccc}1&0&\cdots&0\\{0}&*&\cdots&*\\\vdots&\vdots&{}&\vdots\\{0}&*&\cdots&*\end{array}\right),$$$$U_1=\mbox{matrices of the form}\quad\left(\begin{array}{cc}1&0\\{*}&I_d\end{array}\right),$$$$T_1=\mbox{matrices of the form}\quad\diag(1,*,\ldots,*).$$Then $P_1=L_1U_1$ and $T_1$ is a maximal torus in $L_1\cong{\rm GL}\sb {d}$. We view the character group $X^*(T_1)$ of $T_1$ as the subgroup of $X^*(T)$ generated by $\epsilon_1,\ldots,\epsilon_d$. The morphism ${\rm GL}\sb {d+1}\to\mathbb{A}^{d+1}$, $g\mapsto(g(1,0,\ldots,0))$ induces an isomorphism ${\rm GL}\sb {d+1}/P_1\cong \mathbb{A}^{d+1}$. Over $X^{cone}$ it has the section$$u:X^{cone}\to{\rm GL}\sb {d+1},\quad z\mapsto u(z)=\left(\begin{array}{cc}\Xi_0(z)&\Xi_1(z)\quad\cdots\quad\Xi_{d}(z)\\{0}&I_d\end{array}\right)^{-1}.$$This gives rise to the automorphy factor$$\mu:G\to P_1(H^0(X^{cone},{\mathcal O}_{X^{cone}}))$$$$g\mapsto u(g(z))^{-1}.g.u(z).$$As a matrix valued function on $X^{cone}$ it satisfies the automorphy factor relation$$\mu(gh)(z)=\mu(g)(hz).\mu(h)(z)\quad\mbox{for}\, g,h\in G$$(the three factors are viewed as elements in ${\rm GL}\sb {d+1}(H^0(X^{cone},{\mathcal O}_{X^{cone}}))$, and their product there turns out to ly in $P_1(H^0(X^{cone},{\mathcal O}_{X^{cone}}))$). Let us describe it explicitly. Let $$g=\left(\begin{array}{ccc}a_{00}&\cdots&a_{0d}\\\vdots&{}&\vdots\\a_{d0}&\cdots&a_{dd}\end{array}\right)\in G$$and put\begin{gather}A_j=a_{0j}\Xi_0+\ldots+a_{dj}\Xi_d\quad\mbox{for}\,0\le j\le d;\label{aidef}\end{gather}then\begin{gather}\mu(g^{-1})^{-1}=\left(\begin{array}{cccc}1&0&\cdots&0\\a_{10}A_0^{-1}&a_{11}-a_{10}A_1A_0^{-1}&\cdots&a_{1d}-a_{10}A_dA_0^{-1}\\\vdots&\vdots&{}&\vdots\\a_{d0}A_0^{-1}&a_{d1}-a_{d0}A_1A_0^{-1}&\cdots&a_{dd}-a_{d0}A_dA_0^{-1}\end{array}\right).\label{autexp}\end{gather}In particular we see that the image of $\mu(g^{-1})^{-1}$ under the projection $P_1(H^0(X^{cone},{\mathcal O}_{X^{cone}}))\to L_1(H^0(X^{cone},{\mathcal O}_{X^{cone}}))$ in fact lies in $L_1(H^0(X,{\mathcal O}_{X}))$. We denote it by $\nu(g)$.\\ 

Now let $V\ne0$ be an irreducible $K$-rational representation of $L_1$. For $\mu\in X^*(T_1)$ let $V_{\mu}$ be the maximal subspace of $V$ on which $T_1$ acts through $\mu$. Choose a $\mu=\sum_{i=1}^da_{i}\epsilon_i\in X^*(T_1)$ such that $V_{\mu}\ne 0$ and set $|V|=\sum_{i=1}^da_{i}$; this is independent of the choice of $\mu$, as all $\mu$ with $V_{\mu}\ne 0$ differ by linear combinations of elements of $\Phi$ (see \cite{jan} II.2.2). Viewing $V$ as the constant sheaf with value $V$ on ${{X}}$ we define the coherent ${\mathcal O}_X\otimes_K\widehat{K}$-module$${\mathcal V}_{\widehat{K}}=V\otimes_{{K}}{\mathcal O}_{{X}}\otimes_K\widehat{K}.$$Applying $H^0({X},{\mathcal O}_{{X}}\otimes_K\widehat{K})\otimes_K(.)$ the $K$-linear action of $L_1(K)$ on $V$ gives rise to an action of $L_1(H^0({X},{\mathcal O}_{{X}}\otimes_K\widehat{K}))$ on ${\mathcal V}_{\widehat{K}}$. By the automorphy factor relation we get a $G$-action on ${\mathcal V}_{\widehat{K}}$ by setting\begin{gather}g(f\otimes v)=\widehat{\pi}^{-|V|\omega(\det(g))}f(g^{-1}(.))\nu(g)(1\otimes v)\label{i}\end{gather}for $g\in G$, $v\in V$ and any section $f$ of ${\mathcal O}_{{X}}\otimes_K\widehat{K}$.\\ 

Fix a $L_1/{\mathcal O}_K$-invariant ${\mathcal O}_K$-lattice $V_0$ in $V$ (see \cite{jan} I.10.4). (When we write $V_0$ we always refer to this lattice in $V$; to refer to the weight space in $V$ for the weight $\mu=0$ we reserve the phrase " $V_{\mu}$ for $\mu=0$ ".)  

\begin{lem}\label{intspan} We have $V_0=\oplus_{\mu\in X^*(T_1)}V_{\mu,0}$ with $V_{\mu,0}=V_0\cap V_{\mu}$.
\end{lem}

{\sc Proof:} We reproduce a proof of Schneider and Teitelbaum. Fix $\mu\in X^*(T_1)$. It suffices to construct an element $\Pi_{\mu}$ in the algebra of distributions $\mbox{Dist}(L_1/\mathbb{Z})$ (i.e. defined over $\mathbb{Z}$) which on $V$ acts as a projector onto $V_{\mu}$. For $1\le i\le d$ let $H_i=(de_i)(1)\in\mbox{Lie}(L_1/\mathbb{Z})$; then $d\mu'(H_i)\in\mathbb{Z}$ (inside $\mbox{Lie}(\mathbb{G}_m/\mathbb{Z})$) for any $\mu'\in X^*(T_1)$. According to \cite{hum} Lemma 27.1 we therefore find a polynomial $\Pi\in{\mathbb{Q}}[X_1,\ldots,X_d]$ such that $\Pi(\mathbb{Z}^d)\subset \mathbb{Z}$, $\Pi(d\mu(H_1),\ldots,d\mu(H_d))=1$ and $\Pi(d\mu'(H_1),\ldots,d\mu'(H_d))=0$ for any $\mu'\in X^*(T_1)$ such that $\mu'\ne\mu$ and $V_{\mu'}\ne0$. Moreover \cite{hum} Lemma 26.1 says that $\Pi$ is a $\mathbb{Z}$-linear combination of polynomials of the form $${X_1\choose b_1}\cdots{X_d\choose b_d}\quad\mbox{with integers}\quad b_1,\ldots,b_d\ge0.$$Thus \cite{jan} II.1.12 implies that $$\Pi_{\mu}=\Pi(H_1,\ldots,H_d)$$lies in $\mbox{Dist}(L_1/\mathbb{Z})$. By construction it acts on $V$ as a projector onto $V_{\mu}$.\hfill$\Box$\\  

We denote the pushforward $V\otimes_{{\mathcal O}_K}{\mathcal O}_{\widehat{\mathfrak{X}}}$ of ${\mathcal V}_{\widehat{K}}$ via the specialization map $X\to{\mathfrak{X}}$ again by ${\mathcal V}_{\widehat{K}}$. It is a $G$-equivariant (via formula (\ref{i})) coherent ${\mathcal O}_{\mathfrak{X}}\otimes_{{\mathcal O}_{{K}}}\widehat{K}$-module sheaf on ${\mathfrak{X}}$. 

\begin{satz}\label{intstrex} There exists a $G$-equivariant coherent ${\mathcal O}_{\widehat{\mathfrak{X}}}$-submodule ${\mathcal V}_{{\mathcal O}_{\widehat{K}}}$ of ${\mathcal V}_{\widehat{K}}$ such that ${\mathcal V}_{{\mathcal O}_{\widehat{K}}}\otimes_{{\mathcal O}_{\widehat{K}}}\widehat{K}={\mathcal V}_{\widehat{K}}$. For $Z\in F^0_A$ its restriction to ${\mathfrak{U}}_{Z}$ is \begin{gather} {\mathcal V}_{{\mathcal O}_{\widehat{K}}}|_{{\mathfrak{U}}_{Z}}=\bigoplus_{\mu\in X^*(T_1)}\widehat{\pi}^{(d+1)\overline{\mu}(Z)}V_{\mu,0}\otimes_{{\mathcal O}_K}{\mathcal O}_{\widehat{\mathfrak{U}}_{Z}}.\label{ii}\end{gather} 
\end{satz}

{\sc Proof:} We set\begin{gather}{\mathcal V}_{{\mathcal O}_{\widehat{K}}}|_{\mathfrak Y}=\bigoplus_{\mu\in X^*(T_1)}({\mathcal O}_{\widehat{\mathfrak Y}})^{\overline{\mu}}\otimes_{{\mathcal O}_K}V_{\mu,0}.\end{gather}Clearly $({\mathcal V}_
{{\mathcal O}_{\widehat{K}}}|_{\mathfrak Y})|_{{\mathfrak{U}}_{Z}}$ satisfies formula (\ref{ii}) for any $Z\in F^0_A$, because ${\mathcal J}_{Z}|_{{\mathfrak{U}}_{Z}}=\pi{\mathcal O}_{{\mathfrak{U}}_{Z}}$ and ${\mathcal J}_{Z'}|_{{\mathfrak{U}}_{Z}}={\mathcal O}_{{\mathfrak{U}}_{Z}}$ for $Z'\ne Z$. Since ${\mathfrak X}=G{\mathfrak Y}$ the proof of Theorem \ref{intstrex} is complete once we have Proposition \ref{intgstab} below. \hfill$\Box$\\

\begin{pro}\label{intgstab} Let ${\mathfrak W}_1$, ${\mathfrak W}_2$ be open formal subschemes of ${\mathfrak Y}$, let $g\in G$ such that $g{\mathfrak W}_1={\mathfrak W}_2$. Then the isomorphism$$g:{\mathcal V}_
{{\widehat{K}}}|_{{\mathfrak W}_1}\cong{\mathcal V}_
{{\widehat{K}}}|_{{\mathfrak W}_2}$$induces an isomorphism$$g:({\mathcal V}_
{{\mathcal O}_{\widehat{K}}}|_{\mathfrak Y})|_{{\mathfrak W}_1}\cong({\mathcal V}_
{{\mathcal O}_{\widehat{K}}}|_{\mathfrak Y})|_{{\mathfrak W}_2}.$$
\end{pro}

{\sc Proof:} The crucial arguments for (i)-(iii) below are due to Schneider and Teitelbaum (who considered a similar situation on $X^{cone}$ rather than on $X$ or ${\mathfrak X}$). From Lemma \ref{lbdlext} we deduce: For any open $V\subset{\mathfrak Y}\otimes k$ and any $f\in(({\mathcal O}_{\widehat{\mathfrak Y}})^{\overline{\mu}}\otimes_{{\mathcal O}_{\widehat{K}}}k)(V)$, we have $f=0$ if and only if $f|_{V\cap{\mathfrak{U}}_{Z}}=0$ for any $Z\in F^0_A$. It follows that for any open ${\mathfrak V}\subset{\mathfrak Y}$ and any $f\in(({\mathcal O}_{\widehat{\mathfrak Y}})^{\overline{\mu}}\otimes_{{\mathcal O}_{\widehat{K}}}\widehat{K})({\mathfrak V})$, we have $f\in ({\mathcal O}_{\widehat{\mathfrak Y}})^{\overline{\mu}}({\mathfrak V})$ if and only if $f|_{{\mathfrak V}\cap {\mathfrak{U}}_{Z}}\in({\mathcal O}_{\widehat{\mathfrak Y}})^{\overline{\mu}}({\mathfrak V}\cap {\mathfrak{U}}_{Z})$ for any $Z\in F^0_A$. Since the sum over the $\mu$'s in the definition of ${\mathcal V}_
{{\mathcal O}_{\widehat{K}}}|_{\mathfrak Y}$ is direct, this last statement holds verbatim also with ${\mathcal V}_
{{\mathcal O}_{\widehat{K}}}|_{\mathfrak Y}$ instead of $({\mathcal O}_{\widehat{\mathfrak Y}})^{\overline{\mu}}$.

Assume first that we know Proposition \ref{intgstab} whenever ${\mathfrak W}_1={\mathfrak{U}}_{Z}$, ${\mathfrak W}_2={\mathfrak{U}}_{gZ}$ for some $Z\in F^0_A$ such that $gZ\in F^0_A$. Then consider arbitrary ${\mathfrak W}_1$, ${\mathfrak W}_2$. By construction, ${\mathcal F}_1=g(({\mathcal V}_
{{\mathcal O}_{\widehat{K}}}|_{\mathfrak Y})|_{{\mathfrak W}_1})$ and ${\mathcal F}_2=({\mathcal V}_
{{\mathcal O}_{\widehat{K}}}|_{\mathfrak Y})|_{{\mathfrak W}_2}$ are coherent ${\mathcal O}_{\widehat{\mathfrak W}_2}$-submodules of ${\mathcal V}_
{{\widehat{K}}}|_{{\mathfrak W}_2}$ such that ${\mathcal F}_1\otimes_{{\mathcal O}_{\widehat{K}}}{\widehat{K}}={\mathcal F}_2\otimes_{{\mathcal O}_{\widehat{K}}}{\widehat{K}}={\mathcal V}_
{{\widehat{K}}}|_{{\mathfrak W}_2}$. From our assumption it follows that ${\mathcal F}_1|_{{\mathfrak{U}}_{Z}\cap{\mathfrak W}_2}={\mathcal F}_2|_{{\mathfrak{U}}_{Z}\cap{\mathfrak W}_2}$ for all $Z\in F^0_{A}$. But this indeed implies ${\mathcal F}_1={\mathcal F}_2$ by our remark above.

Now we treat the case ${\mathfrak W}_1={\mathfrak{U}}_{Z}$, ${\mathfrak W}_2={\mathfrak{U}}_{gZ}$ for some $Z\in F^0_{A}$ such that $gZ\in F^0_{A}$. Let us write $Z=t.Y$ with some $t=\diag(t_0,\ldots,t_d)\in T$. Let $x=-\sum_{i=0}^d\omega(t_i)e_i\in A$; this may depend on the choice of $t$, but the group $U_x$ does not, it is canonically associated with $Z$. Also note that $\gamma(x)=\gamma(Z)$ for all $\gamma\in X^*(\overline{T})\otimes{\mathbb{R}}$. Similarly we choose a $gx\in A$ corresponding to $gZ$ in the same manner. (In fact we may view $A$ as an apartment in the {\it extended} building associated with ${\rm GL}_{d+1}(K)$; it is acted on by ${\rm GL}_{d+1}(K)$ and we may take $gx$ to be the image of $x$ under the action of $g$: that this lies again in $A$ follows from our hypothesis $gZ\in F^0_{A}$.)  If $W$ denotes the subgroup of permutation matrices then $N=T\rtimes W$ is the normalizer of $T$ in $G$. By the Bruhat decomposition, there exist $h_x\in U_x$, $h_{gx}\in U_{gx}$ and $n\in N$ such that $g=h_xnh_{gx}$. Therefore we may split up our task into the following cases (i)-(iii):

(i) $g\in T$, (ii) $g\in W$, (iii) $x=gx$ and $g\in U_x$.

(i) Suppose $g=\diag(t_0,\ldots,t_d)$. We claim that in this case $g$ even respects weight spaces; in view of formula (\ref{ii}) this means we must prove that $g$ induces for any $\mu\in X^*(T_1)$ with $V_{\mu}\ne0$ an isomorphism$$g:\widehat{\pi}^{(d+1)\overline{\mu}(x)}V_{\mu,0}\otimes_{{\mathcal O}_K}{\mathcal O}_{\widehat{\mathfrak{U}}_{Z}}\cong\widehat{\pi}^{(d+1)\overline{\mu}(gx)}V_{\mu,0}\otimes_{{\mathcal O}_K}{\mathcal O}_{\widehat{\mathfrak{U}}_{gZ}}.$$From formula (\ref{autexp}) we deduce $\nu(g)=\diag(1,t_1,\ldots,t_d)$. Write $\mu=\sum_{i=1}^da_i\epsilon_i$. Then $\nu(g)$ acts on $V_{\mu}$ as $v\mapsto(\prod_{i=1}^dt_i^{a_i})v$, hence induces an isomorphism $$\nu(g):V_{\mu,0}\cong\pi^kV_{\mu,0}\quad\mbox{with}\quad k=\sum_{i=1}^da_i\omega(t_i).$$But on the other hand $k=-\omega(\mu(g))=|V|\omega(\det(g)^{\frac{1}{d+1}})-\omega(\mu(g))-|V|\omega(\det(g)^{\frac{1}{d+1}})=\omega(\overline{\mu}(g))-|V|\omega(\det(g)^{\frac{1}{d+1}})$ and therefore $$\overline{\mu}(gx)=\overline{\mu}(x)+\omega(\overline{\mu}(g))=\overline{\mu}(Z)+k+|V|\omega(\det(g)^{\frac{1}{d+1}}).$$Together the claim follows by inspecting formula (\ref{i}).

(ii) Now $g\in W$. First consider the case $g(e_0)=e_0$. Then $\mu(g^{-1})^{-1}=g$, hence $\nu(g)$ is a permutation matrix acting on $V$ by isomorphisms $\nu(g):V_{\mu}\cong V_{\nu(g)\mu}$ which restrict to isomorphisms $\nu(g):V_{\mu,0}\cong V_{\nu(g)\mu,0}$ since by assumption $V_0\subset V$ is stable for $L_1({\mathcal O}_K)$. On the other hand $\nu(g)$ acts on $X^*(T_1)$ such that $\mu(x)=(\nu(g)\mu)(gx)$ and hence $\overline{\mu}(x)=\overline{(\nu(g)\mu)}(gx)$ for $\mu\in X^*(T_1)$. It follows that $g$ induces isomorphisms$$\widehat{\pi}^{(d+1)\overline{\mu}(x)}V_{\mu,0}\otimes_{{\mathcal O}_K}{\mathcal O}_{\widehat{\mathfrak{U}}_{Z}}\cong\widehat{\pi}^{(d+1)\overline{(\nu(g)\mu)}(gx)}V_{\nu(g)\mu,0}\otimes_{{\mathcal O}_K}{\mathcal O}_{\widehat{\mathfrak{U}}_{gZ}}$$for any $\mu\in X^*(T_1)$ and we are done for such $g$. Next consider the case where $g\in W$ is the transposition $g=(0i)$ for some $1\le i\le d$. For $0\le j\le d$ let $\Upsilon_j=-\frac{\Xi_j}{\Xi_i}$. Then\begin{gather}\nu(g)=\left(\begin{array}{cc}1&{}\\{}&\begin{array}{ccccccc}1&{}&{}&{}&{}&{}\\{}&\ddots&{}&{}&{}&{}&{}\\{}&{}&1&{}&{}&{}&{}\\\Upsilon_1&\cdots&\Upsilon_{i-1}&\Upsilon_0&\Upsilon_{i+1}&\cdots&\Upsilon_d\\{}&{}&{}&{}&1&{}&{}\\{}&{}&{}&{}&{}&\ddots&{}\\{}&{}&{}&{}&{}&{}&1\end{array}\end{array}\right)\in L_1({\mathcal O}_{\mathfrak X}\otimes_{{\mathcal O}_K}K).\label{bigmat}\end{gather}Let $J=\{1\le j\le d;\,j\ne i\}$. We may factorize $\nu(g)$ as $$\nu(g)=e_i(-\Upsilon_0)\prod_{j\in J}\widetilde{\alpha}_{ij}(\Theta_j)$$with $\Theta_j=\frac{\Xi_j}{\Xi_0}=\frac{\Upsilon_j}{\Upsilon_0}$ and with $\widetilde{\alpha}_{ij}$ as defined by formula (\ref{alphaij}). For the ${\mathcal O}_{\mathfrak X}\otimes_{{\mathcal O}_K}K$-linear extension of the $L_1(K)$-action on $V$ to an $L_1({\mathcal O}_{\mathfrak X}\otimes_{{\mathcal O}_K}K)$-action on ${\mathcal O}_{\mathfrak X}\otimes_{{\mathcal O}_K}V$ we have\begin{gather}\widetilde{\alpha}_{ij}(f)(1\otimes v)\in\sum_{m\ge0}f^m.V_{\mu+m(\epsilon_i-\epsilon_j),0}\label{liearg}\end{gather}for $f\in {\mathcal O}_{\mathfrak X}\otimes_{{\mathcal O}_K}K$, $v\in V_{\mu,0}$ and $j\in J$. To see this define $X_{\alpha}=(d\widetilde{\alpha}_{ij})(1)\in\mbox{Lie}(L_1/\mathbb{Z})$ for $\alpha=\epsilon_i-\epsilon_j$ and then $$X_{\alpha,m}=\frac{X_{\alpha}^m}{m!}\in\mbox{Dist}(L_1/\mathbb{Z})\quad\mbox{for}\quad m\ge0$$(compare \cite{jan} II.1.11 and 1.12). By \cite{jan} II.1.19 we have$$X_{\alpha,m}V_{\mu}\subset V_{\mu+m\alpha}\quad\mbox{and}\quad\widetilde{\alpha}_{ij}(f)(1\otimes v)=\sum_{m\ge0}f^mX_{\alpha,m}(1\otimes v)$$(with $\alpha=\epsilon_i-\epsilon_j$). Since $X_{\alpha,m}$ is defined over $\mathbb{Z}$ we in turn have $X_{\alpha,m}V_{\mu,0}\subset V_{\mu+m\alpha,0}$ and formula (\ref{liearg}) follows. By the above factorization of $\nu(g)$ and formula (\ref{liearg}) we get$$\nu(g)(1\otimes v)\in\sum_{(m_j)_{j\in J}\in\mathbb{N}_0^J}(\mu+\sum_{j\in J}m_j(\epsilon_i-\epsilon_j))(e_i(-\Upsilon_0))(\prod_{j\in J}\Theta_j^{m_j}).V_{\mu+\sum_{j\in J}m_j(\epsilon_i-\epsilon_j),0}$$for $v\in V_{\mu,0}$. Here $\prod_{j\in J}\Theta_j^{m_j}\in\pi^{\sum_{j\in J}m_j(\epsilon_j-\epsilon_0)(gx)}{\mathcal O}_{\mathfrak{X}}({{\mathfrak{U}}_{gZ}})$ and$$(\mu+\sum_{j\in J}m_j(\epsilon_i-\epsilon_j))(e_i(-\Upsilon_0))=(-\Upsilon_0)^{a_i+\sum_{j\in J}m_j}\in\pi^{(a_i+\sum_{j\in J}m_j)(\epsilon_0-\epsilon_i)(gx)}{\mathcal O}_{\mathfrak{X}}({{\mathfrak{U}}_{gZ}})$$if we write $\mu=\sum_{s=1}^da_s\epsilon_s$. Together$$\nu(g)(1\otimes v)\in\sum_{(m_j)_{j\in J}\in\mathbb{N}_0^J}\pi^{\sum_{j\in J}m_j(\epsilon_j-\epsilon_0)(gx)+(a_i+\sum_{j\in J}m_j)(\epsilon_0-\epsilon_i)(gx)}V_{\mu+\sum_{j\in J}m_j(\epsilon_i-\epsilon_j),0}\otimes_{{\mathcal O}_K}{\mathcal O}_{\mathfrak{X}}({{\mathfrak{U}}_{gZ}})$$for $v\in V_{\mu,0}$. We may rewrite the exponent as $\sum_{j\in J}m_j(\epsilon_j-\epsilon_0)(gx)+(a_i+\sum_{j\in J}m_j)(\epsilon_0-\epsilon_i)(gx)=a_i(\epsilon_0-\epsilon_i)(gx)-\sum_{j\in J}m_j(\epsilon_i-\epsilon_j)(gx)$. On the other hand$$\overline{\mu}(x)=\overline{\mu}(gx)-a_i(\epsilon_0-\epsilon_i)(gx)$$$$=\overline{\mu+\sum_{j\in J}m_j(\epsilon_i-\epsilon_j)}(gZ)-a_i(\epsilon_0-\epsilon_i)(gx)+\sum_{j\in J}m_j(\epsilon_i-\epsilon_j)(gx)$$and it follows that $g$ maps the submodule $\widehat{\pi}^{(d+1)\overline{\mu}(x)}V_{\mu,0}\otimes_{{\mathcal O}_K}{\mathcal O}_{\widehat{\mathfrak{U}}_{Z}}$ of ${\mathcal V}_
{{\widehat{K}}}|_{{\mathfrak{U}}_{Z}}$ into $$\bigoplus_{\mu'\in X^*(T_1)}\widehat{\pi}^{(d+1)\overline{\mu'}(gx)}V_{\mu',0}\otimes_{{\mathcal O}_K}{\mathcal O}_{\widehat{\mathfrak{U}}_{gZ}},$$for any $\mu\in X^*(T_1)$. By a symmetry argument (consider $g^{-1}$) we are done for this kind of $g$.

(iii) Now consider the case $x=gx$ and $g\in U_x$. We may assume $g\in U_{\alpha,-\alpha(x)}$ for some $\alpha=\epsilon_i-\epsilon_t\in\Phi$. Thus $g=\widetilde{\alpha}_{it}(u)$ for some $u\in K$ with $\omega(u)\ge-\alpha(x)$. If $i=0$ then $\nu(g)=I_{d+1}$ and our claim is obvious. If $i\ge1$ and $t\ge1$ then $\nu(g)=g$ and it suffices to show that the automorphism $g$ of $V$ induces an automorphism$$g:\bigoplus_{\mu\in X^*(T_1)}\widehat{\pi}^{(d+1)\overline{\mu}(x)}V_{\mu,0}\cong\bigoplus_{\mu\in X^*(T_1)}\widehat{\pi}^{(d+1)\overline{\mu}(x)}V_{\mu,0}.$$But $\omega(u)\ge-\alpha(x)$ implies $\overline{\mu+m\alpha}(x)\le\overline{\mu}(x)+m\omega(u)$ for all $\mu\in X^*(T_1)$, all $m\in\mathbb{N}$ and we conclude as in the proof of formula (\ref{liearg}). Finally assume $t=0$. For $1\le j\le d$ let now $\Upsilon_j=-\frac{u\Xi_j}{\Xi_0+u\Xi_i}$ and $\Upsilon_0=\frac{\Xi_0}{\Xi_0+u\Xi_i}$. Then formula (\ref{bigmat}) holds also in this context. Letting $J=\{1\le j\le d;\,j\ne i\}$ and this time $\Theta_j=-\frac{u\Xi_j}{\Xi_0}=\Upsilon_0^{-1}\Upsilon_j$ we may factorize $\nu(g)$ as $$\nu(g)=e_i(\Upsilon_0)\prod_{j\in J}\widetilde{\alpha}_{ij}(\Theta_j).$$Similarly as in (ii) we get$$\nu(g)(1\otimes v)\in\sum_{(m_j)_{j\in J}\in\mathbb{N}_0^J}(\mu+\sum_{j\in J}m_j(\epsilon_i-\epsilon_j))(e_i(\Upsilon_0))(\prod_{j\in J}\Theta_j^{m_j}).V_{\mu+\sum_{j\in J}m_j(\epsilon_i-\epsilon_j),0}$$for $v\in V_{\mu,0}$ and $\mu\in X^*(T_1)$. Here $\prod_{j\in J}\Theta_j^{m_j}\in\pi^{\sum_{j\in J}m_j(\omega(u)+(\epsilon_j-\epsilon_0)(x))}{\mathcal O}_{\mathfrak{X}}({{\mathfrak{U}}_{Z}})$. On the other hand$$(\mu+\sum_{j\in J}m_j(\epsilon_i-\epsilon_j))(e_i(\Upsilon_0))\in{\mathcal O}_{\mathfrak{X}}({{\mathfrak{U}}_{Z}})$$ because $\Upsilon_0=\frac{\Xi_0}{\Xi_0+u\Xi_i}$ is a unit in ${\mathcal O}_{\mathfrak{X}}({{\mathfrak{U}}_{Z}})$ (because $g=\widetilde{\alpha_{i0}}(u)$ and $g{{\mathfrak{U}}_{Z}}={{\mathfrak{U}}_{Z}}$). Together, observing $\omega(u)\ge-(\epsilon_i-\epsilon_0)(x)$, we obtain$$\nu(g)(1\otimes v)\in\sum_{(m_j)_{j\in J}\in\mathbb{N}_0^J}\pi^{-\sum_{j\in J}m_j(\epsilon_i-\epsilon_j)}V_{\mu+\sum_{j\in J}m_j(\epsilon_i-\epsilon_j),0}\otimes_{{\mathcal O}_K}{\mathcal O}_{\mathfrak{X}}({{\mathfrak{U}}_{Z}})$$for $v\in V_{\mu,0}$. This concludes the proof of Proposition \ref{intgstab} and thus of Theorem \ref{intstrex}.\hfill$\Box$\\

{\it Remark:} If $|V|\in(d+1){\mathbb{Z}}$ (equivalently: if $\overline{\mu}\in X^*({\overline{T}})$ for all $\mu\in X^*(T^1)$ with $V_{\mu}\ne0$) then we could replace the $\widehat{K}$-valued character $\widehat{\pi}^{-|V|\omega(\det(g))}$ in definition (\ref{i}) with the $K$-valued character $\det(g)^{-(d+1)^{-1}|V|}$. Then the scalar extension $K\to\widehat{K}$ could be completely avoided and we would obtain a ${\rm PGL}_{d+1}(K)$-equivariant locally free coherent ${\mathcal O}_{\mathfrak{X}}$-module ${\mathcal V}_{{\mathcal O}_{{K}}}$ such that ${\mathcal V}_{{\mathcal O}_{{K}}}\otimes_{{\mathcal O}_{{K}}}K=V\otimes_{{\mathcal O}_K}{\mathcal O}_{\mathfrak{X}}$.\\

A lattice in $K^{d+1}$ is a free ${\mathcal O}_K$-submodule of $K^{d+1}$ of rank $d+1$. Two lattices $L, L'$ are homothetic if $L'=\lambda L$ for some $\lambda\in K^{\times}$. We denote the homothety class of $L$ by $[L]$. The set of vertices of the Bruhat-Tits building ${\mathcal BT}$ of ${\rm PGL}_{d+1}(K)$ is the set of homothety classes of lattices (always: in $K^{d+1}$). For a lattice chain $$\pi L_s\subsetneq L_1\subsetneq\ldots\subsetneq L_s$$we declare the ordered $s$-tuple $([L_1],\ldots,[L_s])$ to be a pointed $s-1$-simplex (with underlying $s-1$-simplex the unordered set $\{[L_1],\ldots,[L_s]\}$). Call it $\widehat{\eta}$ and consider the set$$N_{\widehat{\eta}}=\{L \mbox{ a lattice in } K^{d+1}\quad|\quad\pi L_s\subsetneq L\subsetneq L_1\}.$$A subset $M_0$ of $N_{\widehat{\eta}}$ is called {\it stable} if for all $L, L'\in M_0$ also $L\cap L'$ lies in $M_0$. We may identify $N_{\widehat{\eta}}$ (and hence $M_0$) with the set of homothety classes which its elements represent; this is independent of the choice of representing lattices for $[L_1],\ldots,[L_s]$. We recall a result from \cite{acy}. Let ${\mathcal F}$ denote a cohomological coefficient system on ${\mathcal BT}$.

\begin{pro}\label{acycri} Let $1\le s\le d$. Suppose that for any pointed $s-1$-simplex $\widehat{\eta}$ with underlying $s-1$-simplex $\eta$ and for any stable subset $M_0$ of $N_{\widehat{\eta}}$ the following subquotient complex of the cochain complex $C^{\bullet}({\mathcal BT},{\mathcal F})$ with values in ${\mathcal F}$ is exact:$${\mathcal F}(\eta)\longrightarrow\prod_{z\in M_0}{\mathcal F}(\{z\}\cup\eta)\longrightarrow\prod_{z,z'\in M_0\atop\{z,z'\}\in F^1}{\mathcal F}(\{z,z'\}\cup\eta).$$Then the $s$-th cohomology group $H^s({\mathcal{BT}},{\mathcal F})$ of $C^{\bullet}({\mathcal BT},{\mathcal F})$ vanishes.\hfill$\Box$\\
\end{pro}
    
A homothety class $[L]$ ( = a vertex of ${\mathcal{BT}}$) corresponds to an irreducible component $Z_{[L]}$ of ${\mathfrak X}\otimes k$. If $L$ represents $[L]$ and $(L/\pi L)^*$ denotes the $k$-vector space dual to $L/\pi L$ then we may regard $Z_{[L]}$ as the successive blowing up of the projective space ${\mathbb P}((L/\pi L)^*)$ in all its $k$-linear subspaces. Now let again $\pi L_s\subsetneq L_1\subsetneq\ldots\subsetneq L_s$ define a pointed $s-1$-simplex $\widehat{\eta}$ and suppose that $L_s={\mathcal O}_K^{d+1}$, the standard lattice. Then $Z_{[L_s]}=Y$, the central component of ${\mathfrak X}\otimes k$ considered in section \ref{cohmodp}. The lattices $L$ with $\pi L_s\subsetneq L\subsetneq L_s$ correspond bijectively to the irreducible components $Z$ of ${\mathfrak X}\otimes k$ with $Z\cap Y\ne\emptyset$, or equivalently to the elements of the set ${\mathcal V}$ from section \ref{cohmodp}. Namely, $L$ defines a $k$-linear subspace in $L_s/\pi L_s$, hence a quotient of $(L_s/\pi L_s)^*$, hence (taking the kernel) a $k$-linear subspace of $(L_s/\pi L_s)^*$, hence a $k$-linear subspace of ${\mathbb P}((L_s/\pi L_s)^*)$; its strict transform under $Y\to {\mathbb P}((L_s/\pi L_s)^*)$ is the element of ${\mathcal V}$ corresponding to $L$. It is clear that a subset $M_0$ of $N_{\widehat{\eta}}$ is stable if and only if the corresponding subset of ${\mathcal V}$ is stable in the sense of section \ref{cohmodp}.

\begin{satz}\label{globnull} Suppose that for every $\mu=\sum_{i=1}^da_i\epsilon_i\in X^*(T_1)$ with $V_{\mu}\ne0$ and for every $1\le j\le d$ we have $\sum_{i\ne j}a_i\le da_j$. Then \begin{gather}H^t(\mathfrak{X},{\mathcal V}_{{\mathcal O}_{\widehat{K}}})=0\label{lifvan}\\H^t(\mathfrak{X},{\mathcal V}_{{\mathcal O}_{\widehat{K}}}\otimes_{{\mathcal O}_{\widehat{K}}}k)=0\label{redvan}\end{gather} for all $t>0$, and \begin{gather}H^0(\mathfrak{X},{\mathcal V}_{{\mathcal O}_{\widehat{K}}})\otimes_{{\mathcal O}_{\widehat{K}}}k=H^0({\mathfrak{X}},{\mathcal V}_{{\mathcal O}_{\widehat{K}}}\otimes_{{\mathcal O}_{\widehat{K}}}k).\label{bac}\end{gather}
\end{satz}

{\sc Proof:} We use the following ordering on $X^*(T)$: Define\begin{gather}\sum_{i=0}^da_i\epsilon_i>\sum_{i=0}^da'_i\epsilon_i\label{ordwei}\end{gather}if and only if there exists a $0\le i_0\le d$ such that $a_i=a_i'$ for all $i<i_0$, and $a_{i_0}>a'_{i_0}$. In particular we get an ordering on $X^*(T_1)$. Then by \cite{jan} II.1.19 the filtration $(F^{\mu}V)_{\mu\in X^*(T_1)}$ of our $K$-rational $L_1$-representation $V$ defined by$$F^{\mu}V=\sum_{\mu'\in X^*(T_1)\atop\mu'\ge\mu}V_{\mu'}$$is stable for the action of $U(K)\cap L_1(K)$. Define the filtration $(F^{\mu}{\mathcal V}_{{\mathcal O}_{\widehat{K}}})_{\mu\in X^*(T_1)}$ of ${\mathcal V}_{{\mathcal O}_{\widehat{K}}}$ by$$F^{\mu}{\mathcal V}_{{\mathcal O}_{\widehat{K}}}={\mathcal V}_{{\mathcal O}_{\widehat{K}}}\cap(F^{\mu}V\otimes_{{\mathcal O}_K}{\mathcal O}_{\widehat{\mathfrak{X}}})$$($\mu\in X^*(T_1)$). Let$$Gr^{\mu}{\mathcal V}_{{\mathcal O}_{\widehat{K}}}=\frac{F^{\mu}{\mathcal V}_{{\mathcal O}_{\widehat{K}}}}{F^{\mu+\epsilon_d}{\mathcal V}_{{\mathcal O}_{\widehat{K}}}}=\frac{F^{\mu}{\mathcal V}_{{\mathcal O}_{\widehat{K}}}}{\bigoplus_{\mu'>\mu}F^{\mu'}{\mathcal V}_{{\mathcal O}_{\widehat{K}}}}.$$Since for any $g\in U(K)$ the automorphy factor $\nu(g)$ is just the image of $g$ under the natural projection $U(K)\to U(K)\cap L_1(K)$ we deduce that the coherent ${\mathcal O}_{\widehat{\mathfrak{X}}}$-modules $F^{\mu}{\mathcal V}_{{\mathcal O}_{\widehat{K}}}$ and $Gr^{\mu}{\mathcal V}_{{\mathcal O}_{\widehat{K}}}$ are $U(K)$-equivariant. Note that the ${U}(K)$-action on $Gr^{\mu}{\mathcal V}_{{\mathcal O}_{\widehat{K}}}$ is with {\it trivial} automorphy factors. We have by construction a canonical isomorphism\begin{gather}Gr^{\mu}{\mathcal V}_{{\mathcal O}_{\widehat{K}}}|_{{\mathfrak Y}}\cong({\mathcal O}_{\widehat{\mathfrak Y}})^{\overline{\mu}}\otimes_{{\mathcal O}_K}V_{\mu,0}.\label{resy}\end{gather}
This is because the composition$$V_{\mu,0}\longrightarrow V_0\cap F^{\mu}V\longrightarrow \frac{V_0\cap F^{\mu}V}{V_0\cap F^{\mu+\epsilon_d}V}$$is an isomorphism.

Now we prove Theorem \ref{globnull}. The hypothesis means that for all $\mu\in X^*(T_1)$ with $V_{\mu}\ne0$ we have $\overline{a}_j(\mu)\le0$ for all $1\le j\le d$. We claim \begin{gather}H^t(\mathfrak{X},Gr^{\mu}{\mathcal V}_{{\mathcal O}_{\widehat{K}}}\otimes_{{\mathcal O}_{\widehat{K}}}k)=0\quad\quad(t>0)\label{grvan}\end{gather}for all $\mu$. First consider the case $\overline{\mu}\in X^*(\overline{T})$. Then Lemma \ref{lbdlext} and equation (\ref{resy}) imply that $Gr^{\mu}{\mathcal V}_{{\mathcal O}_{\widehat{K}}}\otimes_{{\mathcal O}_{\widehat{K}}}k$ is a $U(K)$-equivariant vector bundle on ${\mathfrak X}\otimes k$, of rank ${\rm{rk}}_{{\mathcal O}_K}V_{\mu,0}=\dim_KV_{\mu}$, and that\begin{gather}{\mathcal L}_{Y}(D(\overline{a}(\mu),0,0))\otimes_k(V_{\mu,0}\otimes_{{\mathcal O}_K}k)\cong Gr^{\mu}{\mathcal V}_{{\mathcal O}_{\widehat{K}}}\otimes_{{\mathcal O}_{\widehat{\mathfrak{X}}}}{\mathcal O}_{Y}.\label{denudi}\end{gather}With Corollary \ref{kritvan} (and using $U(K)$-equivariance which tells us that we may assume $W\subset Y$) we get \begin{gather}H^t(\mathfrak{X},Gr^{\mu}{\mathcal V}_{{\mathcal O}_{\widehat{K}}}\otimes_{{\mathcal O}_{\widehat{\mathfrak{X}}}}{\mathcal O}_{W})=0\label{grsimpvan}\end{gather} for all $t>0$, for every component intersection $W\in F^s$ (any $0\le s\le d$). On the other hand, since $Gr^{\mu}{\mathcal V}_{{\mathcal O}_{\widehat{K}}}\otimes_{{\mathcal O}_{\widehat{K}}}k$ is locally free over ${\mathcal O}_{{\mathfrak X}\otimes k}$ we have an exact sequence$$0\longrightarrow Gr^{\mu}{\mathcal V}_{{\mathcal O}_{\widehat{K}}}\otimes_{{\mathcal O}_{\widehat{K}}}k\longrightarrow\prod_{W\in F^0}Gr^{\mu}{\mathcal V}_{{\mathcal O}_{\widehat{K}}}\otimes_{{\mathcal O}_{\widehat{\mathfrak{X}}}}{\mathcal O}_{W}\longrightarrow\prod_{W\in F^1}Gr^{\mu}{\mathcal V}_{{\mathcal O}_{\widehat{K}}}\otimes_{{\mathcal O}_{\widehat{\mathfrak{X}}}}{\mathcal O}_{W}\longrightarrow\ldots.$$Therefore (\ref{grsimpvan}) tells us that to prove (\ref{grvan}) we need to prove $H^s({\mathcal{BT}},{\mathcal F})=0$ for all $s>0$, where ${\mathcal F}$ is the coefficient system on ${\mathcal{BT}}$ defined by$${\mathcal F}(W)=H^0({\mathfrak{X}},Gr^{\mu}{\mathcal V}_{{\mathcal O}_{\widehat{K}}}\otimes_{{\mathcal O}_{\widehat{\mathfrak{X}}}}{\mathcal O}_{W})$$for $W\in F^s$, where we identify $F^s$ with the set of $s$-simplices of ${\mathcal{BT}}$ (any $s$). We apply Proposition \ref{acycri}. Given a pointed $s-1$-simplex $\widehat{\eta}=([L_1],\ldots,[L_s])$ we may, by $G$-equivariance, assume that $[L_s]$ is represented by ${\mathcal O}_K^{d+1}$. Then (\ref{denudi}) and Corollary \ref{kritvan} tell us that the hypothesis of Propositon \ref{acycri} are satisfied, proving (\ref{grvan}). 

The case $\overline{\mu}\notin X^*(\overline{T})$ is easier: there we find\begin{gather}Gr^{\mu}{\mathcal V}_{{\mathcal O}_{\widehat{K}}}\otimes_{{\mathcal O}_{\widehat{K}}}k\cong\prod_{Z\in F^0}\frac{Gr^{\mu}{\mathcal V}_{{\mathcal O}_{\widehat{K}}}\otimes_{{\mathcal O}_{\widehat{\mathfrak{X}}}}{\mathcal O}_{Z}}{{\mathcal O}_{Z}\mbox{-torsion}},\label{zerldelninu}\end{gather}\begin{gather}{\mathcal L}_{Y}(D(\lceil\overline{a}(\mu)\rceil,n(\mu),m(\mu)))\otimes_k(V_{\mu,0}\otimes_{{\mathcal O}_K}k)\cong \frac{Gr^{\mu}{\mathcal V}_{{\mathcal O}_{\widehat{K}}}\otimes_{{\mathcal O}_{\widehat{\mathfrak{X}}}}{\mathcal O}_{Y}}{{\mathcal O}_{Y}\mbox{-torsion}}.\label{deninudi}\end{gather}Combining with Theorem \ref{gallgvan} (and using $U(K)$-equivariance which tells us that it suffices to look at $Z=Y$) we get (\ref{grvan}). 

By the obvious devissage argument we get (\ref{redvan}) from (\ref{grvan}). The base change formula (\ref{bac}) is a consequence of (\ref{lifvan}) for $t=1$, and (\ref{lifvan}) for all $t>0$ is a consequence of (\ref{redvan}) for all $t>0$, as one easily shows by using of the exact sequences$$0\longrightarrow{\mathcal V}_{{\mathcal O}_{\widehat{K}}}\otimes_{{\mathcal O}_{\widehat{K}}}k\stackrel{\pi^r}{\longrightarrow}{\mathcal V}_{{\mathcal O}_{\widehat{K}}}/(\widehat{\pi}^{r+1})\longrightarrow{\mathcal V}_{{\mathcal O}_{\widehat{K}}}/(\widehat{\pi}^{r})\longrightarrow0.$$\hfill$\Box$\\



\begin{thebibliography}{abcdefgh}  

\bibitem{blka}{\it S. Bloch, K. Kato}, $p$-adic \'{e}tale cohomology. Inst. Hautes \'{E}tudes Sci. Publ. Math. No. {\bf 63} (1986), 107--152.
\bibitem{caen}{\it M. Cabanes, M. Enguehard}, Representation theory of finite reductive groups. New Mathematical Monographs, 1. Cambridge University Press, Cambridge (2004)
\bibitem{calu}{\it R. W. Carter, G. Lusztig}, Modular representations of finite groups of Lie type, Proc. London Math. Soc. (3) {\bf 32} (1976), 347--348
\bibitem{dellus}{\it P. Deligne, G. Lusztig}, Representations of reductive groups over finite fields. Ann. of Math. (2) 103 (1976), no. 1, 103--161.
\bibitem{mathan}{\it E. Grosse-Kl\"onne}, Integral structures in automorphic line bundles on the $p$-adic upper half plane, Math. Ann. {\bf 329}, 463--493 (2004)
\bibitem{phien}{\it E. Grosse-Kl\"onne}, Frobenius and Monodromy operators in rigid analysis, and Drinfel'd's symmetric space, Journal of Algebraic Geometry 14 (2005), 391--437
\bibitem{acy}{\it E. Grosse-Kl\"onne}, Acyclic coefficient systems on buildings, to appear in Compositio Math.
\bibitem{latt}{\it E. Grossse-Kl\"onne}, Sheaves of $p$-adic lattices in Weyl modules for ${\rm GL}$, preprint
\bibitem{hart}{\it R. Hartshorne}, Algebraic Geometry, Graduate Texts in Mathematics {\bf 52}, Springer 1977
\bibitem{hum}{\it  J. E. Humphreys}, Introduction to Lie Algebras and Representation Theory. Berlin-Heidelberg-New York: Springer (1972)
\bibitem{ito}{\it T. Ito}, Weight-Monodromy conjecture for $p$-adically uniformized varieties, preprint 2003
\bibitem{jan}{\it J. C. Jantzen}, Representations of algebraic groups, Boston: Academic Press 1987
\bibitem{mus}{\it G. A. Mustafin}, Non-Archimedean uniformization, Math. USSR Sbornik {\bf 34}, 187-214 (1987)
\bibitem{schn}{\it P. Schneider}, The cohomology of local systems on 
$p$-adically uniformized varieties, Math. Ann. {\bf 293}, 623-650 (1992)
\bibitem{schtei}{\it P. Schneider, J. Teitelbaum}, An integral transform for $p$-adic symmetric spaces, Duke Math. J. {\bf 86}, 391-433 (1997)
\bibitem{ast}{\it P. Schneider, J. Teitelbaum}, $p$-adic boundary values. in: Cohomologies $p$-adiques et applications arithmétiques, I. Astérisque No. 278, 51--125 (2002)
\bibitem{teit}{\it J. Teitelbaum}, Modular representations of ${\rm PGL}\sb 2$ and automorphic forms for Shimura curves. Invent. Math. {\bf 113} (1993), no. 3, 561--580.

\end{thebibliography}
\end{document}